\documentclass[11pt,a4paper]{article}
\usepackage{ifpdf}
\pdfoutput=1
\usepackage{amsmath}
\usepackage{amssymb}
\usepackage{mathrsfs}
\usepackage{color}
\usepackage{graphicx}
\usepackage[percent]{overpic}
\usepackage{url}
\usepackage{subfig}
\usepackage{enumerate}
\usepackage{overpic}
\usepackage{amsthm}
\usepackage{moreverb}
\usepackage{authblk}

\def\s{{\bm s}}
\def\b{{\bm b}}

\def\n{{\bm n}}
\def\u{{\bm u}}
\def\v{{\bm v}}

\def\0{\boldsymbol{0}}

\def\mubar{\overline{\mu}}
\def\ubar{\overline{\u}}

\def\qbar{\overline{q}}
\def\dt{\partial_t}

\def\cl {\nonumber \\}
\def\el {\nonumber }

\newtheorem{rem}{Remark}[section]

\newcommand{\bm}[1]{\mbox{\boldmath{$#1$}}}
\def\div{\nabla\cdot}

\usepackage{verbatim}
\usepackage{graphicx}
\usepackage{epstopdf}

\begin{document}
\date{}
\title{Pressure stabilization strategies for a LES filtering Reduced Order Model}


\author[1]{Michele Girfoglio\thanks{mgirfogl@sissa.it}}
\author[2]{Annalisa Quaini\thanks{quaini@math.uh.edu, ORCID ID 0000-0001-9686-9058}}
\author[1]{Gianluigi Rozza\thanks{grozza@sissa.it, ORCID ID 0000-0002-0810-8812}}
\affil[1]{SISSA, International School for Advanced Studies, Mathematics Area, mathLab, via Bonomea, Trieste 265 34136, Italy\\
}
\affil[2]{Department of Mathematics, University of Houston, Houston TX 77204, USA}
\maketitle

\begin{abstract}
We present a stabilized POD-Galerkin reduced order method (ROM) for a Leray model. 
For the implementation of the model, we combine a two-step algorithm called Evolve-Filter (EF) 
with a computationally efficient finite volume method. 
In both steps of the EF algorithm, velocity and pressure fields are approximated using different POD 
basis and coefficients. To achieve pressure stabilization, we consider and compare
 two strategies: the pressure Poisson equation and the supremizer enrichment of the velocity space. 
We show that the evolve and filtered velocity spaces have to be enriched with the supremizer solutions 
related to both evolve and filter pressure fields in order to obtain stable and accurate solutions
with the supremizer enrichment method. 
We test our ROM approach on 2D unsteady flow past a cylinder at Reynolds number $0 \leq Re \leq 100$. 
We find that both stabilization strategies produce comparable errors in the reconstruction
of the lift and drag coefficients, with the pressure Poisson equation method being more
computationally efficient. 


\end{abstract}
\vspace*{0.2cm}
\vspace*{0.5cm}

\section{Introduction}

For about a couple of decades, reduced order models (ROMs) have emerged as an efficient tool for the approximation of
problems governed by parametrized partial differential equations. This success is owed to the fact that
ROMs can significantly reduce the computational 
cost required by classical full order models (FOMs), e.g.~Finite Element Method or Finite Volume Method, 
when several solutions associated to different parameter values are needed.
The basic ROM framework consists of two steps. During a first phase (called \emph{offline}), a database of several solutions is 
collected by solving a FOM of choice for different parameter values. 
Then, during a second phase (called \emph{online}) the database collected in the offline phase is used to compute the solution for 
newly specified values of the parameters in a short amount of time. 
For a comprehensive review on ROMs, we refer the reader to, e.g., 
\cite{hesthaven2015certified, quarteroniRB2016, bennerParSys, Benner2015, Bader2016, ModelOrderReduction}. 

The particular ROM we consider in this paper is based on a POD-Galerkin approach, 
which consists in extracting of the most energetic modes 
representing the system dynamics and projecting the governing equations 
onto the space spanned by these modes. 
For the specific application that we target, i.e. incompressible fluid flow
at moderately high Reynolds number, it is well known that standard POD-Galerkin models lead to instabilities
\cite{Carlberg2013623,Moin1998}. A successful way to cure these instabilities in advection dominated flows
is to adopt subgrid-scale closure models. See, e.g., \cite{wang_turb, Aubry1988}. 
Thus, we choose to work with a Large Eddy Simulation (LES) approach.
 

We focus on a variant of the so-called Leray model \cite{Leray34}, 
where the small-scale effects are described by a set of equations to be added to the discrete 
Navier-Stokes equations. This extra problem acts as a differential low-pass filter \cite{abigail_CMAME}.
For its  actual implementation, we use the Evolve-Filter (EF) algorithm 
\cite{Boyd1998283, Fischer2001265, Dunca2005,layton_CMAME}.
One of the novelties of our approach is that we use a Finite Volume (FV) method \cite{Girfoglio2019,Girfoglio2021},
while the vast majority of the works on Leray-type models use a Finite Element framework.
In the context of regularized ROMs, the Leray model and the EF algorithm have been 
thoroughly investigated. See, e.g., \cite{Gunzburger2019b,Xie2018,Wells2017,Gunzburger2019,Xie2016}. 
In all of these works, the filtering approach is only employed at reduced order level.
In \cite{Girfoglio2020}, for the first time the LES filtering is used also at the FOM level, i.e.~to generate the snapshot data.
Such an approach provides a ROM which is fully consistent with the FOM as 
the same mathematical framework is used during both the \emph{offline} and \emph{online} stages. 

This paper could be seen as an extension of \cite{Girfoglio2020}. Therein, we used a
pressure Poisson equation (PPE) in the online stage \cite{akhtar2009stability,Stabile2018}
as a pressure stabilization technique. Here, we compare the PPE method with 
the supremizer enrichment of the velocity space \cite{Rozza2007,Ballarin2014,Stabile2018,Gerner2011},
i.e.~another techniques that provides pressure stability.
The main objective of this work is to test the accuracy and efficiency of these two methods within our LES filtering approach. 
We show that adapting the supremizer enrichment to the EF algorithm
is not a trivial exercise. Indeed, the supremizer method becomes accurate only 
when the evolve and filtered velocity spaces are enriched 
by the supremizer solutions associated to both evolve and filter pressure fields.    
Since we are interested in the reconstruction of the pressure field at reduced order level,
time is the only parameter we consider. We vary no physical and/or geometrical parameters. 
We test our framework on 2D flow past a cylinder with time-dependent Reynolds number $0 \leq Re(t) \leq 100$ \cite{turek1996, John2004}. 

The work is organized as follows. Sec.~\ref{sec:FOM} describes 
the full order model and the numerical method we use for it. In
Sec.~\ref{sec:ROM}, we describe the reduced order model. 
The numerical examples are reported in Sec. \ref{sec:results}. Finally,
Sec. \ref{sec:conclusions} provides conclusions and perspectives.

\section{The full order model}\label{sec:FOM}

In this section, we describe our Full Order Model (FOM).  
We consider a fixed domain $\Omega \subset \mathbb{R}^D$ with $D = 2, 3$ over a time 
interval of interest ($t_0$, $T$) $\subset \mathbb{R}^+$. 
The so-called \emph{Leray model} couples the Navier-Stokes equations (NSE) with a differential filter as follows:
\begin{align}
& \rho\, \dt \u + \rho\,\div \left(\ubar \otimes \u\right) - 2\mu \Delta\u + \nabla p = 0 &\text{in }\Omega \times (t_0,T), \label{eq:leray1}  \\
& \div \u = 0&\text{in }\Omega \times (t_0,T),  \label{eq:leray2}  \\
& -2 \alpha^2\Delta\ubar +\ubar +\nabla \lambda = \u&\text{in }\Omega \times (t_0,T), \label{eq:leray3}  \\
& \div \ubar = 0&\text{in }\Omega \times (t_0,T), \label{eq:leray4}
\end{align}
where $\rho$ is the fluid density, $\mu$ is the dynamic viscosity, $\u$ is the fluid velocity, 
$p$ is the fluid pressure, $\ubar$ is the \emph{filtered} velocity, and
variable $\lambda$ is a Lagrange multiplier to enforce the incompressibility constraint for $\ubar$. 
Parameter $\alpha$ in eq.~\eqref{eq:leray3} can be interpreted as a \emph{filtering radius}.
In this paper, $\alpha$ will be constant in space and time. More sophisticated choices
are possible but will not be considered here. 
Problem \eqref{eq:leray1}-\eqref{eq:leray4} is endowed with 
suitable boundary conditions
\begin{align}
&\u = \u_D  &&\mbox{on } \partial\Omega_D \times(t_0,T), \label{eq:bc-d} \\ 
&(2\mu \nabla\u - p\mathbf{I})\n = \0  &&\mbox{on } \partial\Omega_N \times(t_0,T), \label{eq:bc-n} \\ 
&\ubar = \u_D  &&\mbox{on } \partial\Omega_D \times(t_0,T), \label{eq:bc-filter-d} \\ 
&(2\alpha^2 \nabla\ubar - \lambda\mathbf{I})\n = \0  &&\mbox{on } \partial\Omega_N \times(t_0,T).\label{eq:bc-filter-n}
\end{align}
and the initial data $\u= \u_0$ in $\Omega \times\{t_0\}$. Here $\overline{\partial\Omega_D}\cup\overline{\partial\Omega_N}=\overline{\partial\Omega}$ and $\partial\Omega_D \cap\partial\Omega_N=\emptyset$. In addition, $\u_D$ and $\u_0$ are given.
Note that we restrict our attention to homogeneous Neumann boundary conditions. Of course, the methodology we propose 
can be extended to non-homogeneous Neumann conditions, as well as other boundary conditions (e.g., Robin conditions).
 \cite{BQV} 

\subsection{The Evolve-Filter algorithm}\label{sec:FOM_2}

We start with the time discretization of model \eqref{eq:leray1}-\eqref{eq:leray4}.
Let $\Delta t \in \mathbb{R}$, $t^n = t_0 + n \Delta t$, with $n = 0, ..., N_T$ and $T = t_0 + N_T \Delta t$. Moreover, we denote by $y^n$ the approximation of a generic quantity $y$ at the time $t^n$. 

To decouple the Navier-Stokes system \eqref{eq:leray1}-\eqref{eq:leray2} from the filter system \eqref{eq:leray3}-\eqref{eq:leray4} 
at each time step, we consider the Evolve-Filter (EF) algorithm \cite{Boyd1998283, Fischer2001265, Dunca2005}. This algorithm 
reads as follows: given $\u^{n-1}$ and $\u^{n}$, at time $t^{n+1}$:
\begin{itemize}
\item[-] \textit{Evolve}: find intermediate velocity and pressure $(\v^{n+1},q^{n+1})$ such that
\begin{align}
&\rho\, \frac{3}{2\Delta t}\, \v^{n+1} + \rho\, \div \left(\u^* \otimes \v^{n+1}\right) - 2\mu\Delta\v^{n+1} +\nabla q^{n+1} = \b^{n+1},\label{eq:evolve-1.1}\\
& \div \v^{n+1} = 0\label{eq:evolve-1.2},
\end{align}
with boundary conditions
\begin{align}
&\v^{n+1} = \u_D^{n+1}  \quad \mbox{on } \partial\Omega_D \times(t_0,T), \label{eq:bc-d_1} \\ 
&(2\mu \nabla\v^{n+1} - q^{n+1}\mathbf{I})\n = \0  \quad\mbox{on } \partial\Omega_N \times(t_0,T), \label{eq:bc-n_1} 
\end{align}
and initial condition $\v^0 = \u_0$ in $\Omega \times\{t_0\}$. In eq.~\eqref{eq:evolve-1.1}, we set 
$\u^* = 2 \u^n-\u^{n-1}$ and $\b^{n+1} = (4\u^n - \u^{n-1})/(2\Delta t)$. In \eqref{eq:evolve-1.2},
we have used Backward Differentiation Formula of order 2 (BDF2) \cite{quarteroni2007numerical} for
the discretization of the Eulerian time derivative.
\item[-] \textit{Filter}: find $(\u^{n+1},\lambda^{n+1})$ such that
\begin{align}
&-\alpha^2 \Delta\u^{n+1} +\u^{n+1} +\nabla \lambda^{n+1} = \v^{n+1}, \label{eq:evolve-2.1}\\
& \div \u^{n+1}  = 0 \label{eq:filter-1.2},
\end{align}
with boundary conditions
\begin{align}
& \u^{n+1} = \u_D^{n+1}  \quad \mbox{on } \partial\Omega_D \times(t_0,T), \label{eq:bc-filter-d_1} \\ 
& (2\alpha^2 \nabla\u^{n+1} - \lambda^{n+1}\mathbf{I})\n = \0  \quad\mbox{on } \partial\Omega_N \times(t_0,T).\label{eq:bc-filter-n_1}
\end{align}
\end{itemize}
We consider $\u^{n+1}$ and $q^{n+1}$ the approximation of the fluid velocity and pressure at the time $t^{n+1}$, respectively.


\begin{rem}\label{rem:gen_Stokes}
Filter problem \eqref{eq:evolve-2.1}-\eqref{eq:filter-1.2} can be considered 
a generalized Stokes problem. In fact, if we multiply eq.~\eqref{eq:evolve-2.1}
by $\rho/\Delta t$ and rearrange the terms we obtain:
\begin{align}
\frac{\rho}{\Delta t} \u^{n+1}  - \mubar \Delta\u^{n+1} + \nabla \qbar^{n+1} & = \frac{\rho}{\Delta t} \v^{n+1}, \quad \mubar = \rho \frac{\alpha^2}{\Delta t}, \label{eq:filter-1.1}
\end{align}
where $\qbar^{n+1} = \rho \lambda^{n+1}/\Delta t$ can be seen as a filtered pressure. 
Problem \eqref{eq:filter-1.1},\eqref{eq:filter-1.2} can be
seen as a time dependent Stokes problem discretized by the Backward Euler (or BDF1) scheme
and with viscosity $\mubar$. 
A solver for problem \eqref{eq:filter-1.1},\eqref{eq:filter-1.2} can be obtained 
by adapting a standard linearized Navier-Stokes solver. 
\end{rem}



\subsection{Space discrete problem by a Finite Volume method}\label{sec:FOM_3}
For the space discretization of problems \eqref{eq:evolve-1.1}-\eqref{eq:evolve-1.2} and \eqref{eq:filter-1.1},\eqref{eq:filter-1.2}, 
we adopt a FV method. We partition the computational domain $\Omega$ into cells or control volumes $\Omega_i$, with $i = 1, \dots, N_{c}$, where $N_{c}$ is the total number of cells in the mesh. 
Let  \textbf{A}$_j$ be the surface vector of each face of the control volume, 
with $j = 1, \dots, M$. 

The fully discretized form of problem \eqref{eq:evolve-1.1}-\eqref{eq:evolve-1.2} is given by
\begin{align}
&\rho\, \frac{3}{2\Delta t}\, \v^{n+1}_i + \rho\, \sum_j^{} \varphi^*_j \v^{n+1}_{i,j} - 2\mu \sum_j^{} (\nabla\v^{n+1}_i)_j \cdot \textbf{A}_j + \sum_j^{} q^{n+1}_{i,j} \textbf{A}_j  = {\bm b}^{n+1}_i, \label{eq:disc_evolve1} \\
&\sum_j^{} (\nabla q^{n+1})_j \cdot \textbf{A}_j = \sum_j^{} (\textbf{H}(\v_i^{n+1}))_j \cdot \textbf{A}_j, \label{eq:disc_evolve2}
\end{align}
with $\varphi^*_j = \u^{*}_j \cdot \textbf{A}_j$ and
\begin{align}
\textbf{H}(\v^{n+1}_i) = -\rho \sum_j^{} \varphi^*_j \v^{n+1}_{i,j} + 2\mu \sum_j^{} (\nabla\v^{n+1}_i)_j \cdot \textbf{A}_j + {\bm b}^{n+1}_i. \label{eq:H}
\end{align}
In \eqref{eq:disc_evolve1}-\eqref{eq:H}, $\v^{n+1}_i$ and ${\bm b}^{n+1}_i$ denote the average velocity and source term in control volume $\Omega_i$, respectively. Moreover, we denote with $\v^{n+1}_{i,j}$ and $q^{n+1}_{i,j}$ the velocity and pressure
associated to the centroid of face $j$ normalized by the volume of $\Omega_i$.

The fully discrete problem associated to the filter problem \eqref{eq:filter-1.1},\eqref{eq:filter-1.2} is given by 
\begin{align}
&\frac{\rho}{\Delta t} \u^{n+1}_i - \mubar\sum_j^{} (\nabla\u^{n+1}_i)_j \cdot \textbf{A}_j +  \sum_j^{} \qbar^{n+1}_{i,j} \textbf{A}_j = \frac{\rho}{\Delta t} \v^{n+1}_i, \label{eq:disc_filter1} \\
&\sum_j^{} (\nabla \qbar^{n+1}_i)_j \cdot \textbf{A}_j = \sum_j^{} (\overline{\textbf{H}}(\u^{n+1}_i))_j \cdot \textbf{A}_j, \label{eq:disc_filter2}
\end{align}
\begin{align}
\overline{\textbf{H}}(\u^{n+1}_i) =  \mubar\sum_j^{}(\nabla\u^{n+1}_i)_j \cdot \textbf{A}_j + \dfrac{\rho}{\Delta t}\v^{n+1}_i. \label{eq:Hbar}
\end{align}
In \eqref{eq:disc_filter1}-\eqref{eq:Hbar}, we denoted with $\u^{n+1}_i$ the average filtered velocity 
in control volume $\Omega_i$, while $\qbar^{n+1}_{i,j}$ is the auxiliary pressure
at the centroid of face $j$ normalized by the volume of $\Omega_i$.

For more details on the full discretization of both problems \eqref{eq:disc_evolve1}-\eqref{eq:H}
and \eqref{eq:disc_filter1}-\eqref{eq:Hbar}, we refer the reader to \cite{Girfoglio2019}.

We have implemented the EF algorithm within the C++ finite volume library OpenFOAM\textsuperscript{\textregistered} \cite{Weller1998}.
For the solution of the linear system associated with \eqref{eq:disc_evolve1}-\eqref{eq:disc_evolve2} we used the PISO algorithm \cite{PISO}, 
while for problem \eqref{eq:disc_filter1}-\eqref{eq:disc_filter2} we chose a slightly modified version of the SIMPLE algorithm
\cite{SIMPLE}, called SIMPLEC algorithm \cite{Doormaal1984}. Both PISO and SIMPLEC are partitioned algorithms
that decouple the computation of the pressure from the computation of the velocity.

\section{The reduced order model}\label{sec:ROM}
The Reduced Order Model (ROM) we propose is based on a framework introduced in \cite{Girfoglio2020}. 
Here in Sec~\ref{sec:ROM_1} we 
describe the procedure we use to construct a POD-Galerkin ROM, while in Sec.~\ref{sec:ROM_2} 
we present two different strategies for pressure stabilization at reduced order level. 
Finally, Sec.~\ref{sec:ROM_3} describes the method we apply to
enforce non-homogeneous Dirichlet boundary conditions  \eqref{eq:bc-d}, \eqref{eq:bc-filter-d}
at the reduced order level. 
The ROM computations are carried out using ITHACA-FV \cite{RoSta17}, 
an in-house open source C++ library. 

\subsection{A POD-Galerkin projection method}\label{sec:ROM_1}


We approximate velocity fields $\v$ and $\u$ and pressure fields $q$ and $\overline{q}$ 
as linear combinations of the dominant modes (basis functions), which are assumed to be dependent on space variables only,
multiplied by scalar coefficients that depend on the time: 
\begin{align}
\v \approx \v_r = \sum_{i=1}^{N_{v_r}} \beta_i(t) \bm{\varphi}_i(\bm{x}), \quad 
q \approx q_r = \sum_{i=1}^{N_{q_r}} \gamma_i(t) \psi_i(\bm{x}), \label{eq:ROM_1} \\
\u \approx \u_r = \sum_{i=1}^{N_{u_r}} \overline{\beta_i}(t) \bm{\overline{\varphi}}_i(\bm{x}), \quad
\overline{q} \approx \overline{q}_r = \sum_{i=1}^{N_{{\overline{q}_r}}} \overline{\gamma_i}(t) \overline{\psi}_i(\bm{x}). \label{eq:ROM_2}
\end{align}
In \eqref{eq:ROM_1}-\eqref{eq:ROM_2}, $N_{\Phi_r}$ denotes the cardinality of a reduced basis for the space field $\Phi$ belongs to.
Note that ${\Phi}$ could be either a scalar or a vector field. 
Using \eqref{eq:ROM_1} to approximate $\v^{n+1}$ and $q^{n+1}$ in \eqref{eq:evolve-1.1}-\eqref{eq:evolve-1.2},
we obtain
\begin{align}
\rho\, \frac{3}{2\Delta t}\, \v_r^{n+1} + \rho\, \div \left(\u_r^* \otimes \v_r^{n+1}\right) - 2\mu\Delta\v_r^{n+1} +\nabla q_r^{n+1} = \b_r^{n+1}, \label{eq:red-1.1} \\
 \div \v_r^{n+1} = 0, \label{eq:red-1.2}
\end{align}
where $\u_r^* = 2 \u_r^n-\u_r^{n-1}$ and $\b_r^{n+1} = (4\u_r^n - \u_r^{n-1})/(2\Delta t)$.
Then, using \eqref{eq:ROM_2} to approximate $\u^{n+1}$ and $\overline{q}^{n+1}$ in \eqref{eq:filter-1.1},\eqref{eq:filter-1.2}
we get:
\begin{align}
\frac{\rho}{\Delta t} \u_r^{n+1}  - \mubar \Delta\u_r^{n+1} + \nabla \qbar_r^{n+1} = \frac{\rho}{\Delta t} \v_r^{n+1}, \label{eq:red-2.1} \\
 \div \u_r^{n+1} = 0. \label{eq:red-2.2} 
\end{align}

Several techniques to generate the reduced basis spaces are available in the literature, e.g.~the
Proper Orthogonal Decomposition (POD), the Proper Generalized Decomposition (PGD) 
and the Reduced Basis (RB) with a greedy sampling strategy.
See, e.g., \cite{Rozza2008, ChinestaEnc2017, Kalashnikova_ROMcomprohtua, quarteroniRB2016, Chinesta2011, Dumon20111387, Huerta2020, ModelOrderReduction}. 
We find the reduced basis by using the method of snapshots.
To this purpose, we solve the FOM described in Sec.~\ref{sec:FOM}  
for each time $t^k \in \{t^1, \dots, t^{N_s}\} \subset (t_0, T]$. 
The snapshots matrices are obtained from the full-order snapshots:  
\begin{align}\label{eq:space}
\bm{\mathcal{S}}_{{{\Phi}}} = [{{\Phi}}(t^1), \dots, {{\Phi}}(t^{N_s})] \in \mathbb{R}^{N_{\Phi_h} \times N_s} \quad
\text{for} \quad {{\Phi}} = \{\v, \u, q, \overline{q}\},
\end{align}
where the subscript $h$ denotes a solution computed with the FOM and $N_{\Phi_h}$
is the dimension of the space field $\Phi$ belong to in the FOM. 
The POD problem consists in finding, for each value of the dimension of the POD space $N_{POD} = 1, \dots, N_s$, the scalar coefficients $a_1^1, \dots, a_1^{N_s}, \dots, a_{N_s}^1, \dots, a_{N_s}^{N_s}$ and functions ${\bm{\zeta}}_1, \dots, {\bm{\zeta}}_{N_s}$ 
that minimize the error between the snapshots and their projection onto the POD basis. In the $L^2$-norm, we have
\begin{align}
E_{N_{POD}} = \text{arg min} \sum_{i=1}^{N_s} ||{{\Phi}_i} - \sum_{k=1}^{N_{POD}} a_i^k {\bm{\zeta}}_k || \quad \forall N_{POD} = 1, \dots, N_s    \cl
\text{with} \quad ({\bm{\zeta}}_i, {\bm{\zeta}}_j)_{L_2(\Omega)} = \delta_{i,j} \quad \forall i,j = 1, \dots, N_s. \label{eq:min_prob}
\end{align}

It can be shown \cite{Kunisch2002492} that eq.~\eqref{eq:min_prob} is equivalent to the following eigenvalue problem
\begin{align}
\bm{\mathcal{C}}^{{\Phi}} \bm{Q}^{{\Phi}} &= \bm{Q}^{{\Phi}} \bm{\Lambda}^{{\Phi}}, \label{eq:eigen_prob} \\
\mathcal{C}_{ij}^\Phi &= ({\Phi}(t^i), {\Phi}(t^j))_{L_2(\Omega)} \quad \text{for} \quad i,j = 1, \dots, N_s,
\end{align}
where $\bm{\mathcal{C}}^{{\Phi}}$ is the correlation matrix computed from the snapshot matrix $\bm{\mathcal{S}}_{{{\Phi}}}$, $\bm{Q}^{{\Phi}}$ is the matrix of eigenvectors and $\bm{\Lambda}^{{\Phi}}$ is a diagonal matrix whose diagonal entries are the
eigenvalues of $\bm{\mathcal{C}}^{{\Phi}}$. 
Then, the basis functions are obtained as follows:
\begin{align}\label{eq:basis_func}
{\bm{\zeta}}_i = \dfrac{1}{N_s \Lambda_i^\Phi} \sum_{j=1}^{N_s} {\Phi}_j Q_{ij}^\Phi.
\end{align}
The POD modes resulting from the aforementioned methodology are:
\begin{align}\label{eq:spaces}
L_\Phi = [\bm{\zeta}_1, \dots, \bm{\zeta}_{N_{\Phi_r}}] \in \mathbb{R}^{N_{\Phi_h} \times N_{\Phi_r}},
\end{align}
where $N_{\Phi_r} < N_s$ are chosen according to the eigenvalue decay. 
The reduced order model can be obtained through a Galerkin projection of the governing equations onto the POD spaces. 

Let
\begin{align}
&M_{r_{ij}} = (\bm{\varphi}_i, \bm{\varphi}_j)_{L_2(\Omega)}, \quad \widetilde{M}_{r_{ij}} = (\bm{\varphi}_i, \overline{\bm{\varphi}_j})_{L_2(\Omega)}, \quad A_{r_{ij}} = (\bm{\varphi}_i, \Delta \bm{\varphi}_j)_{L_2(\Omega)}, \label{eq:matrices_evolve1} \\
&B_{r_{ij}} = (\bm{\varphi}_i, \nabla \psi_j)_{L_2(\Omega)}, \quad P_{r_{ij}} = (\psi_i, \nabla \cdot \bm{\varphi}_j)_{L_2(\Omega)},
\label{eq:matrices_evolve2}
\end{align}
where $\bm{\varphi}_i$ and $\psi_i$ are the basis functions in \eqref{eq:ROM_1}. 
At time $t^{n+1}$, the reduced algebraic system 
for problem \eqref{eq:red-1.1}-\eqref{eq:red-1.2}
is: 
\begin{align}
&\rho\, \frac{3}{2\Delta t} \bm{M}_r \bm{\beta}^{n+1} + \rho \bm{G}_r(\overline{\bm{\beta}}^{n}, \overline{\bm{\beta}}^{n-1}) \bm{\beta}^{n+1} - 2\mu \bm{A}_r \bm{\beta}^{n+1} + \bm{B}_r \bm{\gamma}^{n+1} = \dfrac{\rho}{\Delta t} \widetilde{\bm{M}}_r \left(2\overline{\bm{\beta}}^{n} - \dfrac{1}{2}\overline{\bm{\beta}}^{n-1}\right), \label{eq:reduced_1} \\
&\bm{P}_r \bm{\beta}^{n+1} = 0, \label{eq:reduced_2}
\end{align}
where vectors $\bm{\beta}^{n+1}$ and $\bm{\gamma}^{n+1}$ contain the values of coefficients $\beta_i$ and $\gamma_i$ 
in \eqref{eq:ROM_1} at time $t^{n+1}$.
The term $\bm{G_r}(\overline{\bm{\beta}}^{n}, \overline{\bm{\beta}}^{n-1}) \bm{\beta}^{n+1}$ in \eqref{eq:reduced_1}
is related to the non-linear convective term: 
\begin{align}\label{eq:convective_matrix}
\left(\bm{G_r}(\overline{\bm{\beta}}^{n}, \overline{\bm{\beta}}^{n-1}) \bm{\beta}^{n+1}\right)_i = (2\overline{\bm{\beta}}^{n} - \overline{\bm{\beta}}^{n-1})^T  \bm{\mathcal{G}}_{r_{i\bm{..}}} \bm{\beta}^{n+1}
\end{align}
where $\bm{\mathcal{G}_r}$ is a third-order tensor defined as follows 
\begin{align}\label{eq:C_tensor}
\mathcal{G}_{r_{ijk}} = (\bm{\varphi_i}, \nabla \cdot (\bm{\varphi_j} \otimes \overline{\bm{\varphi}_k}))_{L_2(\Omega)}.
\end{align}
See \cite{quarteroni2007numerical, Rozza2009} for more details.

Next, let 
\begin{align}
&{\overline{M}}_{r_{ij}} = (\bm{\overline{\varphi}}_i, \bm{\overline{\varphi}}_j)_{L_2(\Omega)}, \quad
\overline{A}_{r_{ij}} = (\bm{\overline{\varphi}}_i, \Delta \bm{\overline{\varphi}}_j)_{L_2(\Omega)}, \label{eq:matrices_filter1} \\
&\overline{B}_{r_{ij}} = (\bm{\overline{\varphi}}_i, \nabla \overline{\psi}_j)_{L_2(\Omega)}, \quad
\overline{P}_{r_{ij}} = (\overline{\psi}_i, \nabla \cdot \bm{\overline{\varphi}}_j)_{L_2(\Omega)},    \label{eq:matrices_filter2}
\end{align}
where $\overline{\bm{\varphi}}_i$ and $\overline{\psi}_i$ are the basis functions in \eqref{eq:ROM_2}. 
At time $t^{n+1}$, the reduced algebraic system 
for problem \eqref{eq:red-2.1}-\eqref{eq:red-2.2} is
\begin{align}
&\frac{\rho}{\Delta t} \bm{{\overline{M}}}_r \bm{\overline{\beta}}^{n+1}  - \mubar \bm{\overline{A}}_r \bm{\overline{\beta}}^{n+1} + \bm{\overline{B}}_r\bm{\overline{\gamma}}^{n+1} = \frac{\rho}{\Delta t} \bm{\widetilde{M}}^T_r \bm{\beta}^{n+1} \label{eq:reduced2_1}, \\
&\bm{\overline{P}}_r \bm{\overline{\beta}}^{n+1} = 0. \label{eq:reduced2_2}
\end{align}
where vectors $\overline{\bm{\beta}}^{n+1}$ and $\overline{\bm{\gamma}}^{n+1}$ contain the values of coefficients $\overline{\beta}_i$ and $\overline{\gamma}_i$ in \eqref{eq:ROM_2} at time $t^{n+1}$.

Finally, the initial conditions for the ROM algebraic system \eqref{eq:reduced_1}-\eqref{eq:reduced_2}
are obtained performing a Galerkin projection of the initial full order condition onto the POD basis spaces:
\begin{align}
{\beta^0}_i = (\v(\bm{x},t_0), \bm{\varphi}_i)_{L_2(\Omega)}, \quad
{\overline{\beta}^0}_i = (\u(\bm{x},t_0), \bm{\overline{\varphi}}_i)_{L_2(\Omega)}. \el
\end{align}

The complete reduced algebraic system at time $t^{n+1}$ is given by \eqref{eq:reduced_1}-\eqref{eq:reduced_2},\eqref{eq:reduced2_1}-\eqref{eq:reduced2_2}.


\subsection{Pressure fields reconstruction and pressure stability}\label{sec:ROM_2}
It is well known that the reduced problem \eqref{eq:reduced_1}, \eqref{eq:reduced_2}, \eqref{eq:reduced2_1}, \eqref{eq:reduced2_2} presents stability issues because the approximation spaces need to satisfy the inf-sup (Ladyzhenskaya-Brezzi-Babuska) condition \cite{BREZZI199027, boffi_mixed}.
In a standard finite element (FE) NSE framework, the inf-sup condition reads: 
\begin{equation}\label{eq:inf-sup}
\inf_{q_h \in \mathcal{Q}_h} \sup_{\u_h \in \mathcal{V}_h} \dfrac{<\div \u_h, q_h>}{||\div \u_h|| ||q_h||} \geq \gamma > 0,
\end{equation}
where $\mathcal{Q}_h$ is the FE space for the pressure approximation,
$\mathcal{V}_h$ is the FE space for the approximation of the velocity field,
and $\gamma$ is a constant that does not depend on the mesh size $h$.
 In order to obtain a stable and accurate reconstruction of the pressure field at the reduced level, different approaches have been proposed. 
One option is to use a global POD basis for both pressure and velocity field and same temporal coefficients
\cite{Bergmann2009Enablers,Lorenzi2016}. Another option is represented by the supremizer enrichment method; see, e.g., \cite{Rozza2007,Ballarin2014,Stabile2018,Gerner2011}.
Finally, one can take the divergence of the momentum equation to obtain a Poisson equation for the pressure 
that is projected onto a POD basis; see, e.g.,
\cite{akhtar2009stability,Stabile2018}. 
This third method is called Poisson pressure equation (PPE).

In this work, we test and compare the performances of two methods: 
the PPE method (already combined with the EF algorithm in \cite{Girfoglio2020})
and the supremizer enrichment method (not yet tested for the EF algorithm).
As we will show, the extension of the supremizer enrichment method to the 
EF algorithm is not straightforward. 

\subsubsection{Pressure Poisson Equation method}

We take the divergence of eq.~\eqref{eq:evolve-1.1} and account for conditions~\eqref{eq:evolve-1.2} 
to obtain the Poisson pressure equation for the Evolve step:
\begin{align}
&\Delta q^{n+1} = -\rho\, \nabla \cdot \left(\div \left(\u^* \otimes \v^{n+1}\right)\right). \label{eq:system_def1_2}
\end{align}
So, at the Evolve step instead of solving \eqref{eq:evolve-1.1}-\eqref{eq:evolve-1.2}, 
we solve the modified systems \eqref{eq:evolve-1.1}, \eqref{eq:system_def1_2} with boundary condition
\eqref{eq:bc-d_1} and 
\begin{align}
& \partial_{{n}} q^{n+1} = -2 \mu \bm{n} \cdot \left(\nabla \times \nabla \times \v^{n+1} \right) - \bm{n} \cdot \left(\rho\dfrac{3}{2\Delta t}\v^{n+1} - \bm{b}^{n+1}\right) \quad \mbox{on } \partial\Omega_N \times(t_0,T), \label{eq:system_def2_2}
\end{align}
where $\partial_{{n}}$ denotes the derivative with respect to outgoing normal \bm{n}. 

We proceed similarly for the Filter step. We take the divergence of eq.~\eqref{eq:filter-1.1} and account for condition \eqref{eq:filter-1.2}
to obtain the Poisson pressure equation:
\begin{align}
&\Delta \overline{q}^{n+1} = 0. \label{eq:system_def1_4}
\end{align}
The Filter step becomes solving \eqref{eq:filter-1.1}, \eqref{eq:system_def1_4} with 
boundary condition \eqref{eq:bc-filter-d_1} and
\begin{align}
& \partial_{{n}} \overline{q}^{n+1} = -2 \overline{\mu} \bm{n} \cdot \left(\nabla \times \nabla \times \u^{n+1} \right)  \quad \mbox{on } \partial\Omega_N \times(t_0,T), \label{eq:system_def2_4}
\end{align}

The reader interested in enforcing non-homogeneous Neumann conditions for the pressure fields
is referred to \cite{Orszag1986, JOHNSTON2004221}.

\begin{rem}\label{rem2}
Systems \eqref{eq:evolve-1.1}-\eqref{eq:evolve-1.2} and \eqref{eq:filter-1.1},\eqref{eq:filter-1.2}
are not equivalent to systems \eqref{eq:evolve-1.1},\eqref{eq:system_def1_2} and \eqref{eq:filter-1.1},\eqref{eq:system_def1_4}
for steady flows \cite{Li2020, Orszag1986, JOHNSTON2004221}. As discussed in Remark \ref{rem:gen_Stokes}, 
the filter problem can be seen as a time-dependent Stokes problem. 
\end{rem}

%
%

Using \eqref{eq:ROM_1} to approximate $\v^{n+1}$ and $q^{n+1}$ in \eqref{eq:system_def1_2},
we obtain
\begin{align}
&\Delta q^{n+1}_r = -\rho\, \nabla \cdot \left(\div \left(\u^*_r \otimes \v^{n+1}_r\right)\right). \label{eq:system_def1_2r}
\end{align}
After space discretization, the matrix form of eq.~\eqref{eq:system_def1_2r} reads:
\begin{align}
&\bm{D}_r \bm{\gamma}^{n+1} + \rho \bm{J}_r(\overline{\bm{\beta}}^{n}, \overline{\bm{\beta}}^{n-1}) \bm{\beta}^{n+1} - 2\mu \bm{N}_r\bm{\beta}^{n+1} - \dfrac{\rho}{2\Delta t}\left(3\bm{F}_r\bm{\beta}^{n+1} - 4 \overline{\bm{F}}_r\bm{\beta}^{n} + \overline{\bm{F}}_r\bm{\beta}^{n-1}\right) = 0 \label{eq:reduced_q},
\end{align}
where 
\begin{align}
&D_{r_{ij}} = (\nabla \psi_i, \nabla \psi_j)_{L_2(\Omega)}, \quad N_{r_{ij}} = (\bm{n} \times \nabla \psi_i, \nabla \times \bm{\varphi}_j) _{L_2(\partial \Omega)}, \label{eq:pressure_matrices1} \\
& F_{r_{ij}} = (\psi_i, \bm{n} \cdot \bm{\varphi}_j)_{L_2(\partial\Omega)}, \quad \overline{F}_{r_{ij}} = (\psi_i, \bm{n} \cdot \overline{\bm{\varphi}}_j)_{L_2(\partial\Omega)}. \label{eq:pressure_matrices2} 
\end{align}
The residual associated with the non-linear term in the equation \eqref{eq:reduced_q} is evaluated with the same strategy
used for eq.~\eqref{eq:reduced_1}. We have
\begin{align}\label{eq:convective_matrix_2}
\left(\bm{J_r}(\overline{\bm{\beta}}^{n}, \overline{\bm{\beta}}^{n-1}) \bm{\beta}^{n+1}\right)_i = (2\overline{\bm{\beta}}^{n} - \overline{\bm{\beta}}^{n-1})^T  \bm{\mathcal{J}}_{r_{i\bm{..}}} \bm{\beta}^{n+1},
\end{align}
where $\bm{\mathcal{J}}_{r}$ is a third-order tensor defined as follows
\begin{align}\label{eq:J_tensor}
\mathcal{J}_{r_{ijk}} = (\nabla \psi_i, \nabla \cdot (\bm{\varphi}_j \otimes \overline{\bm{\varphi}}_k))_{L_2(\Omega)}.
\end{align}

Using \eqref{eq:ROM_2} to approximate $\overline{q}^{n+1}$ in
eq.~\eqref{eq:system_def1_4}, we get 
\begin{align}
&\Delta \overline{q}^{n+1}_r = 0. \label{eq:system_def1_4r}
\end{align}
Once discretized in space, eq.~\eqref{eq:system_def1_4r} can be written in  the matrix form as:
\begin{align}
&\overline{\bm{D}}_r \overline{\bm{\gamma}}^{n+1} - 2\overline{\mu} \overline{\bm{N}}_r \overline{\bm{\beta}}^{n+1} = 0, \label{eq:reduced_q_bar}
\end{align}
where
\begin{align}
& \overline{D}_{r_{ij}} =(\nabla \overline{\psi}_i, \nabla \overline{\psi}_j)_{L_2(\Omega)},  \quad \overline{N}_{r_{ij}} = (\bm{n} \times \nabla \overline{\psi}_i, \nabla \times \bm{\overline{\varphi}_j}) _{L_2(\partial \Omega)}. \label{eq:pressure_matrices3}
\end{align}

With the PPE method, at every time step the ROM algebraic system that has to be solved is \eqref{eq:reduced_1}, \eqref{eq:reduced_q}, \eqref{eq:reduced2_1}, \eqref{eq:reduced_q_bar}.


\subsubsection{Supremizer enrichment method}
For a given pressure basis function, the supremizer is the solution that permits the realization of the inf-sup condition \eqref{eq:inf-sup}. 
One has to find the supremizer for each pressure basis function.
Here, we use an approximated supremizer enrichment procedure:
instead of pressure basis functions we use pressure snapshots. 
This procedure allows to drastically reduce the \emph{online} computational cost. 
In fact, the supremizer basis functions do not depend on the particular pressure basis functions 
but are computed directly from the pressure snapshots during the \emph{offline} phase. 
The downside of this approximated procedure is that it is not possible to rigorously show that the inf-sup condition is satisfied.
One only relies on heuristic criteria or checks during a post-processing stage.

For NSE, the supremizer $\s(t^i)$ relative to pressure snapshot $q(t^i)$
 is found by solving the following problem: 
\begin{align}
&\Delta {\s}(t^i) = -\nabla q(t^i) && \mbox{in } \Omega, \label{eq:sup1}\\
&{\s}(t^i) = 0 &&\mbox{on } \partial\Omega,  \label{eq:sup2}
\end{align}
for $i = 1, \dots, N_s$. 
For more details, the reader is referred to \cite{Rozza2007,Ballarin2014,Stabile2018,Gerner2011}.
For the EF algorithm, in addition to \eqref{eq:sup1}-\eqref{eq:sup2} 
one has to solve an analogous problem  for each auxiliary pressure snapshot $\overline{q}(t^i)$:
\begin{align} 
&\Delta \overline{\s}(t^i) = -\nabla \overline{q}(t^i) && \mbox{in } \Omega,  \label{eq:sup3}\\
&\overline{\s}(t^i) = 0 &&\mbox{on } \partial\Omega,  \label{eq:sup4}
\end{align}
for $i = 1, \dots, N_s$.
Once problems \eqref{eq:sup1}-\eqref{eq:sup2} and \eqref{eq:sup3}-\eqref{eq:sup4} are solved, two snapshots matrices of supremizer solutions are assembled: 
\begin{align}
\bm{\mathcal{S}}_{{{\s}}} = [{{\s}}(t^1), \dots, {{\s}}(t^{N_s})] \in \mathbb{R}^{N_{\v_h} \times N_s}, \quad \bm{\mathcal{S}}_{{{\overline{\s}}}} = [{{\overline{\s}}}(t^1), \dots, {{\overline{\s}}}(t^{N_s})] \in \mathbb{R}^{N_{\u_h} \times N_s}. \el
\end{align}
A POD procedure is applied to the matrices above in order to obtain the supremizer POD basis functions: 
\begin{align}
L_\s = [\bm{\chi}_1, \dots, \bm{\chi}_{N_{\s_r}}] \in \mathbb{R}^{N_{\v_h} \times N_{\s_r}}, \quad
L_{\overline{\s}} = [\bm{\overline{\chi}}_1, \dots, \bm{\overline{\chi}}_{N_{{\overline{\s}_r}}}] \in \mathbb{R}^{N_{\u_h} \times N_{{\overline{\s}_r}}}, \el 
\end{align}
where $N_{\s_r}$ and $N_{{\overline{\s}_r}}$ conform to the notation introduced in Sec.~\ref{sec:ROM_1}.
These velocity supremizer basis functions are added to the reduced velocity spaces, which become:
\begin{align}
\tilde{L}_\v = [\bm{\varphi}_1, \dots, \bm{\varphi}_{N_{\v_r}}, \bm{\chi}_1, \dots, \bm{\chi}_{N_{\s_r}}] \in \mathbb{R}^{N_{\v_h} \times (N_{\v_r} + N_{\s_r})}, \label{eq:sup_POD3}\\
\tilde{L}_\u = [\bm{\overline{\varphi}}_1, \dots, \bm{\overline{\varphi}}_{N_{\u_r}}, \bm{\overline{\chi}}_1, \dots, \bm{\overline{\chi}}_{N_{{\overline{\s}_r}}}] \in \mathbb{R}^{N_{\u_h} \times (N_{\u_r} + N_{{\overline{\s}_r}})}, \label{eq:sup_POD4}
\end{align}

Let us call $\text{SUP}_1$ the supremizer enrichment method described above. In Sec.~\ref{sec:results}, we will show that
$\text{SUP}_1$ does not lead to an accurate 
reconstruction of the pressure fields. To increase the accuracy, 
we add the supremizer solutions associated to both pressure fields
to the evolve and filtered velocity spaces, i.e.:
\begin{align}
\tilde{L}_\v = [\bm{\varphi}_1, \dots, \bm{\varphi}_{N_{\v_r}}, \bm{\chi}_1, \dots, \bm{\chi}_{N_{\s_r}}, \bm{\overline{\chi}}_1, \dots, \bm{\overline{\chi}}_{N_{{\overline{\s}_r}}}] \in \mathbb{R}^{N_{\v_h} \times (N_{\v_r} + N_{\s_r} + N_{{\overline{\s}_r}})}, \label{eq:sup_POD5}\\
\tilde{L}_\u = [\bm{\overline{\varphi}}_1, \dots, \bm{\overline{\varphi}}_{N_{\u_r}}, \bm{\overline{\chi}}_1, \dots, \bm{\overline{\chi}}_{N_{{\overline{\s}_r}}},  \bm{\chi}_1, \dots, \bm{\chi}_{N_{\s_r}}] \in \mathbb{R}^{N_{\u_h} \times (N_{\u_r} + N_{{\overline{\s}_r}} + N_{\s_r})}. \label{eq:sup_POD6}
\end{align}
We call $\text{SUP}_2$ this realization of the supremizer enrichment.

With either $\text{SUP}_1$ or $\text{SUP}_2$, the ROM algebraic system that has to be solved at every time step is \eqref{eq:reduced_1}, \eqref{eq:reduced_2}, \eqref{eq:reduced2_1}, \eqref{eq:reduced2_2}.

\subsection{Treatment of the Dirichlet boundary conditions: the lifting function method}\label{sec:ROM_3}

In order to homogeneize the velocity fields snapshots and make them independent of the boundary conditions,
we use the lifting function method \cite{Stabile2018}. The lifting functions  
are problem-dependent: they have to be divergence free in order to retain 
the divergence-free property of the basis functions and they have to satisfy the boundary conditions of the FOM. 

The velocity snapshots are modified as follows:
\begin{align}
\v'_h = \v_h - \sum_{j=1}^{N_{BC}} v_{{BC}_j}(t)\bm{\chi}_j(\bm{x}), \quad
\u'_h = \u_h - \sum_{j=1}^{N_{BC}} u_{{BC}_j}(t)\bm{\chi}_j(\bm{x}), \el
\end{align}
where $N_{BC}$ is the number of non-homogeneous Dirichlet boundary conditions, $\bm{\chi}_j (\bm{x})$ 
are the lifting functions, and $v_{{BC}_j}$ and $u_{{BC}_j}$ are suitable temporal coefficients. 
The POD is applied to the snapshots satisfying the  homogeneous boundary conditions. 
Then, the boundary value is added back in this way:
\begin{align}
\v_r= \sum_{j=1}^{N_{BC}} v_{{BC}_j}(t)\bm{\chi}_j(\bm{x})  + \sum_{i=1}^{N_{v_r}} \beta_i(t) \bm{\varphi}_i(\bm{x}), \quad
\u_r= \sum_{j=1}^{N_{BC}} u_{{BC}_j}(t)\bm{\chi}_j(\bm{x})  + \sum_{i=1}^{N_{u_r}} \overline{\beta_i}(t) \overline{\bm{\varphi}_i}(\bm{x}). \el
\end{align}

\section{Numerical results}\label{sec:results}

We consider 2D flow past a cylinder \cite{John2004,turek1996}. 
This well-know benchmark will allow us to assess the
ability of the ROM approaches presented in Sec.~\ref{sec:ROM}
to reconstruct the time evolution of the velocity and pressure fields. 
We have thoroughly investigated the 2D flow past a cylinder with a finite volume FOM in \cite{Girfoglio2019}
and with a ROM using a PPE approach in \cite{Girfoglio2020}. 
Here, we aim at comparing the supremizer enrichment method and the PPE method 
in terms of 
errors and speed-up.

The computational domain is a 2.2 $\times$ 0.41 rectangular channel with a cylinder of radius 0.05 centered at (0.2, 0.2), 
when taking the bottom left corner of the channel as the origin of the axes. 
The channel is filled with fluid with density $\rho = 1$ and viscosity $\mu = 10^{-3}$.
We impose a no slip boundary condition on the upper and lower walls and on the cylinder. At the inflow, we prescribe
the following velocity profile:
\begin{align}\label{eq:cyl_bc}
\v(0,y,t) = \left(\dfrac{6}{0.41^2} \sin\left(\pi t/8 \right) y \left(0.41 - y \right), 0\right), \quad y \in [0, 0.41], \quad t \in (0, 8],
\end{align}
and ${\partial q}/{\partial \n} = {\partial \overline{q}}/{\partial \n} = 0$. At the outflow, we prescribe $\nabla \v \cdot  \n = 0$ and $q = \qbar = 0$. 
We start the simulations from fluid at rest. Note that the Reynolds number is time dependent, with $0 \leq Re(t) \leq 100$ \cite{turek1996}.

\begin{figure}[htb!]
\centering
\includegraphics[height=0.115\textwidth]{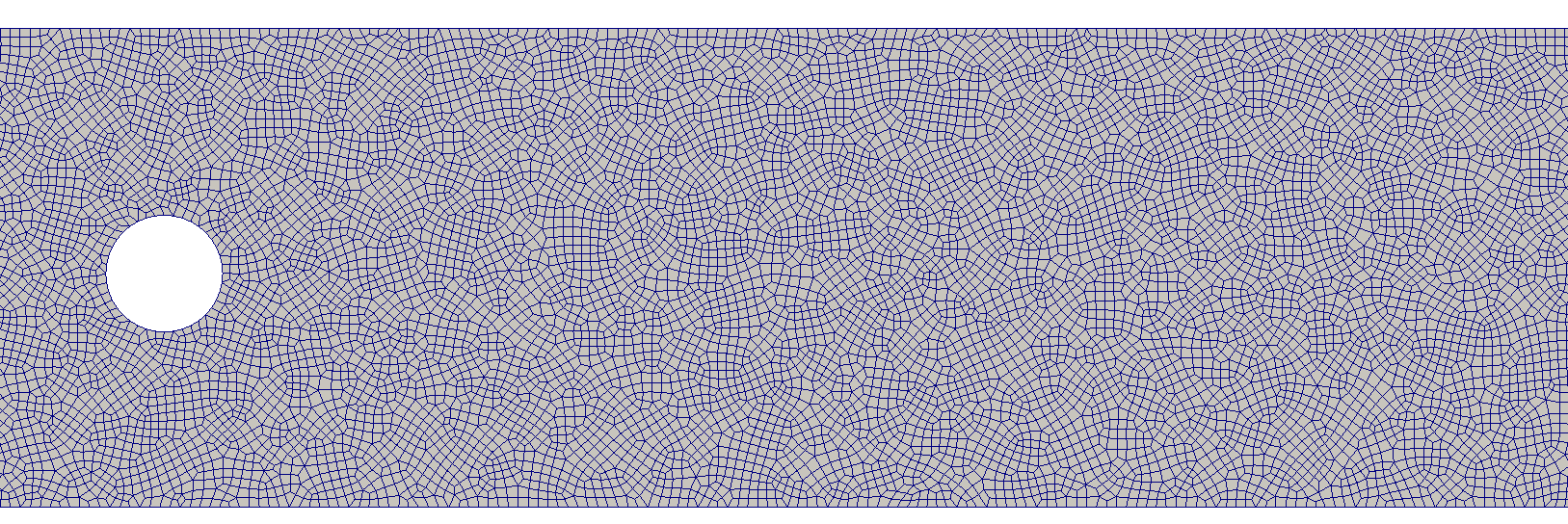}
\includegraphics[height=0.11\textwidth]{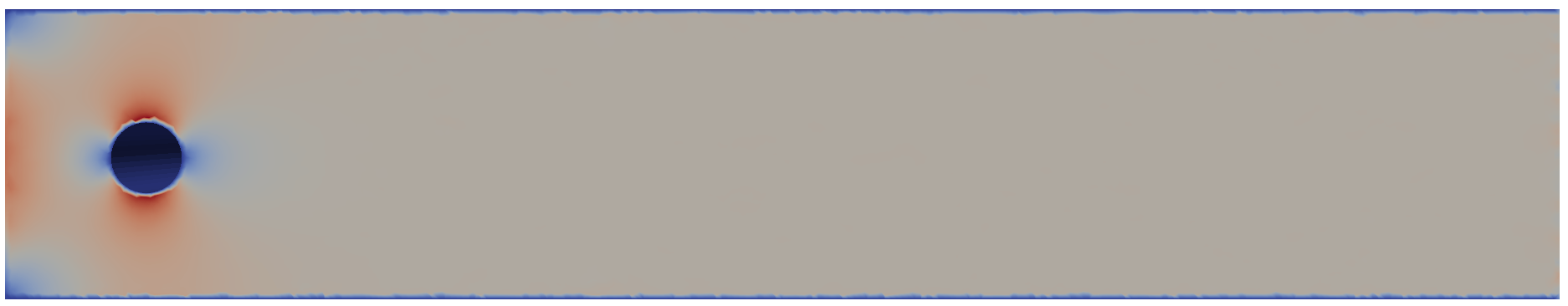}
\includegraphics[height=0.11\textwidth]{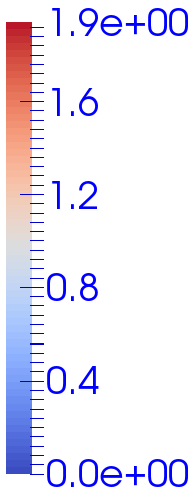}
\caption{Part of the mesh under consideration (left) and lifting function for velocity (right).}
\label{fig:example_cyl}
\end{figure}

We consider a hexaedral computational grid with $h_{min} =  4.2e-3$ and $h_{max} = 1.1e-2$
for a total of $1.59e4$ cells. The quality of this mesh is high: it features
very low values of maximum non-orthogonality (36$^\circ$), 
average non-orthogonality (4$^\circ$), skewness (0.7), and maximum aspect ratio (2). 
Fig.~\ref{fig:example_cyl} (left) shows a part of the mesh. 
We chose this mesh because it is the coarsest among all the meshes considered in \cite{Girfoglio2019}
and thus the most challenging for our LES approach. 
We fix the time step to $\Delta t = 4e-4$ for both FOM and ROM simulations. 
For the convective term, we use a second-order accurate Central Differencing scheme \cite{Lax1960}. In this
way, we avoid introducing stabilization and are thus able to assess the effect of the
filter. 

We need to enforce a time-dependent  non-uniform Dirichlet boundary condition at the inlet. 
For this purpose, we consider a divergence free function with the following non-uniform velocity distribution:
\begin{equation}\label{eq:BC_lift}
\bm{\chi} (0, y) = \left(\dfrac{6}{0.41^2} y \left(0.41 - y \right), 0\right), \quad y \in [0, 2.2],
\end{equation}
and uniform null values on the rest of the boundary. See Figure \ref{fig:example_cyl} (right).

We will compare our findings with those in \cite{Stabile2018}. The choice of Ref.~\cite{Stabile2018}
is due to the fact that therein the authors develop
a NSE-ROM finite volume framework both with PPE and supremizer enrichment methods 
for the reconstruction of the pressure field. 

We set $\alpha = 0.0032$ and refer to \cite{Girfoglio2019} for details on this choice. 
The snapshots are collected every 0.1 s (i.e., $N_s = 80$) using an equispaced grid method in time. 
Therefore, the dimension of the correlation matrix 
$\bm{\mathcal{C}}^{{\Phi}}$ in \eqref{eq:eigen_prob} is $80 \times 80$. 
For a convergence test as the number of snapshots increases, see \cite{Girfoglio2020}. 
Table \ref{tab:cum2D} reports the first 4 cumulative eigenvalues for the velocity, pressure, and supremizer fields. 
In order to retain 99\% of the energy for the ROM, we need 2 modes for $\v$, 2 modes for $q$, 
2 modes for $\u$, and 1 mode for $\overline{q}$. 
As for $\text{SUP}_1$, we consider a number of supremizer modes greater than pressure modes as suggested 
in \cite{Stabile2018}: 4 modes for $\s$ and 3 modes for $\overline{\s}$. On the other hand, for $\text{SUP}_2$ we take into account an equal number of pressure and supremizer modes: 2 modes for $\s$ and 1 mode for $\overline{\s}$.

 \begin{table}[htp!]
\centering
\begin{tabular}{lcccccc}
\multicolumn{2}{c}{} \\
\cline{1-7}
N modes & $\u$ & $\v$ & $q$ & $\overline{q}$ & $\s$ & $\overline{\s}$ \\
\hline
 1 & 0.999588 & 0.999582 & 0.967431 & 0.999985 & 0.736795 & 0.999899\\
 2 & 0.999924 & 0.999924 & 0.999916 & 0.999997 & 0.999594 & 0.999986\\
 3 & 0.999998 & 0.999998 & 0.999995 & 0.999999 & 0.999977 & 0.999999\\
 4 & 0.999999 & 0.999999 & 0.999998 & 0.999995  & 0.999988 & 0.999999\\
\hline
\end{tabular}
\caption{First 4 cumulative eigenvalues for the velocity, pressure and supremizer fields.}
\label{tab:cum2D}
\end{table}

\begin{figure}[htb!]
\centering
 \begin{overpic}[width=0.45\textwidth]{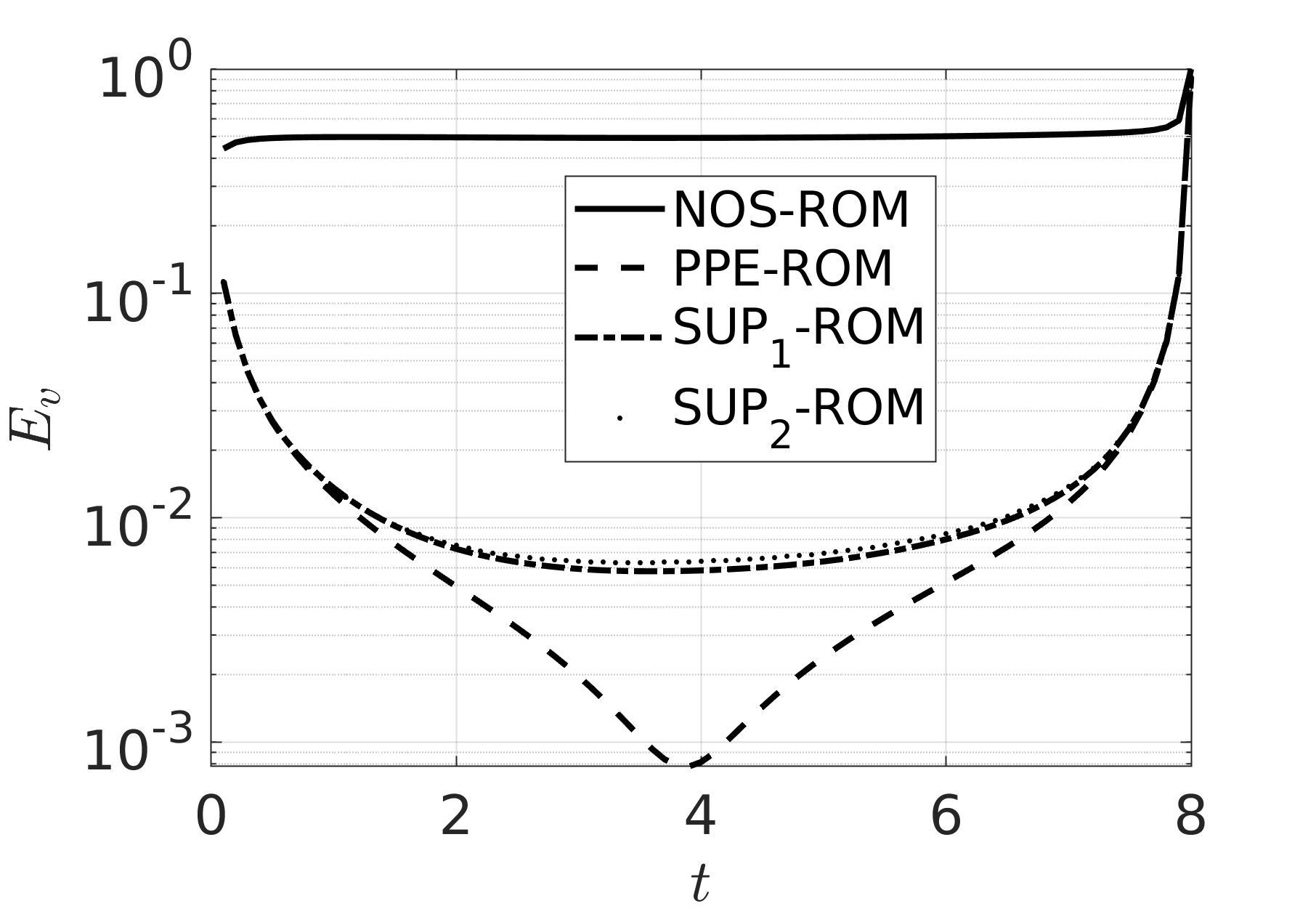}
      \end{overpic}
\begin{overpic}[width=0.45\textwidth]{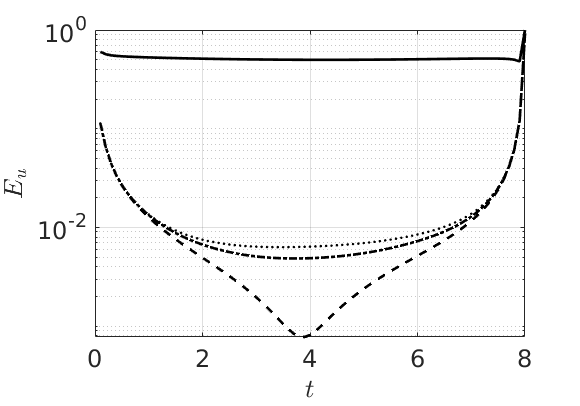}
      \end{overpic}\\
\begin{overpic}[width=0.45\textwidth]{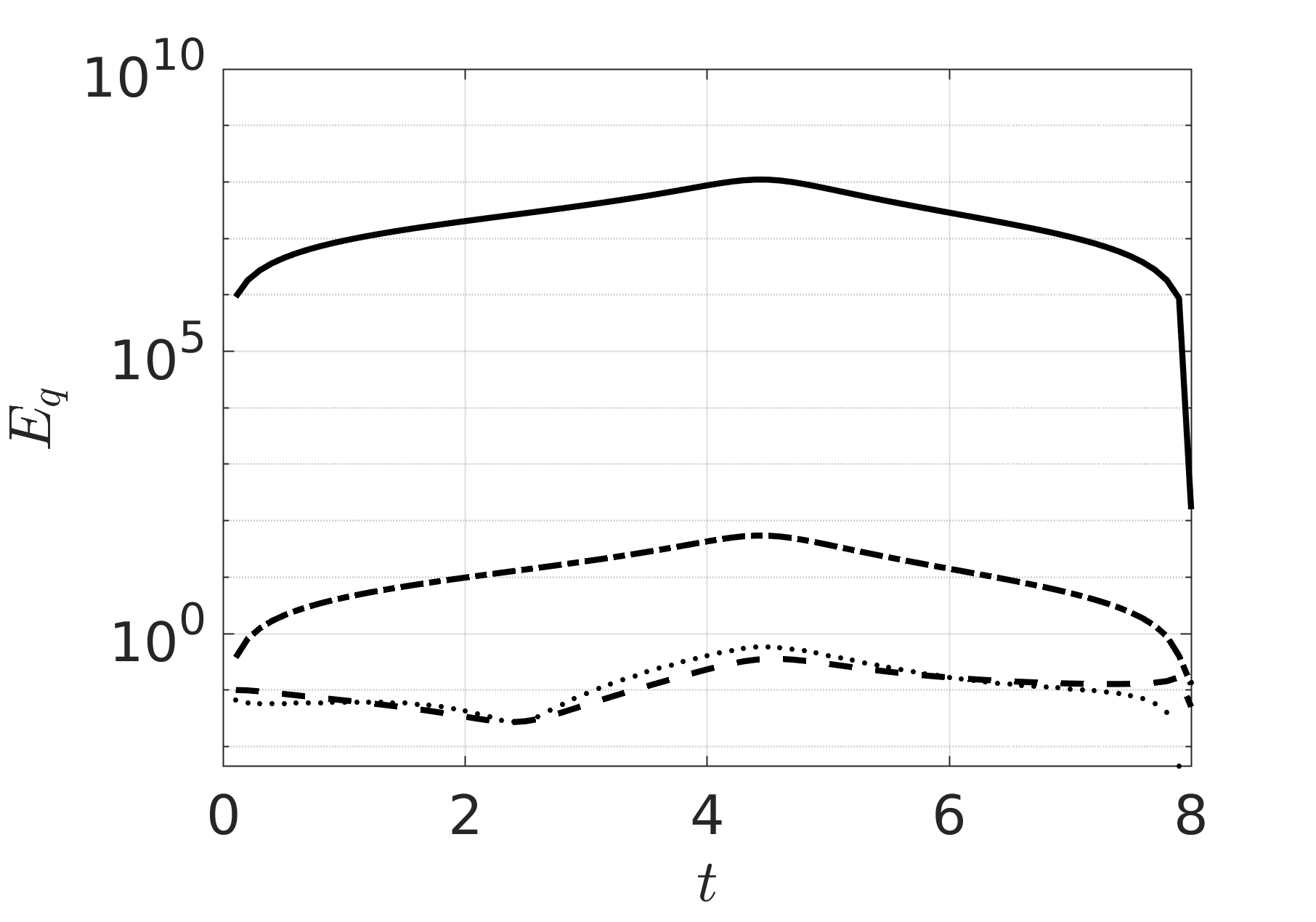}
      \end{overpic}
\begin{overpic}[width=0.45\textwidth]{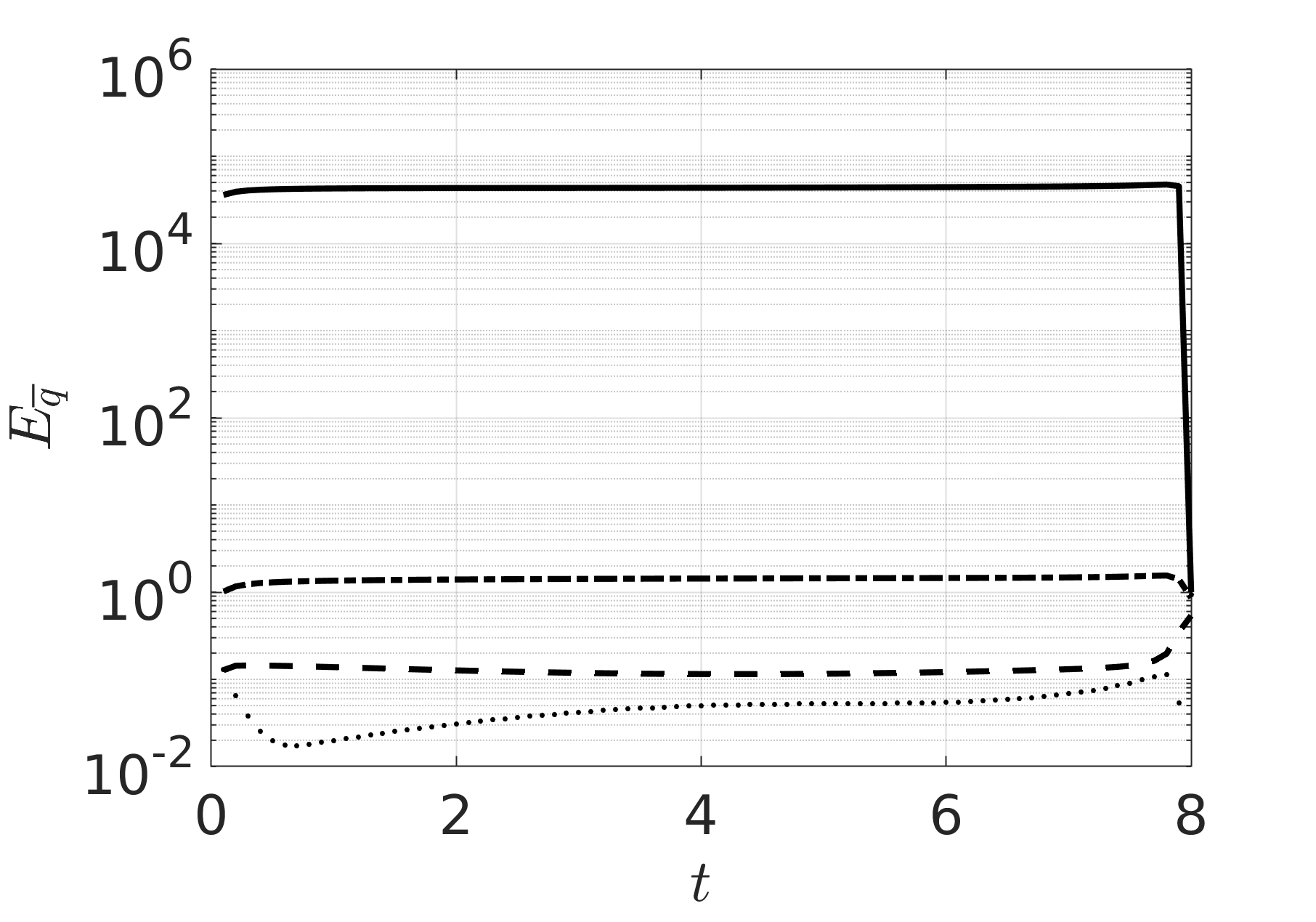}
      \end{overpic}\\
\caption{Time history of $L^2$ norm of the relative error for the velocity fields (top) 
and pressure fields (bottom) for NOS-ROM, PPE-ROM, $\text{SUP}_1$-ROM, and $\text{SUP}_2$-ROM. 
The legend in the top left panel si common to all the panels. 
}
\label{fig:err_t}
\end{figure}

We calculate the $L^2$ relative error: 
\begin{equation}\label{eq:error1}
E_{\Phi} = \dfrac{||\Phi_h(t) - \Phi_r(t)||_{L^2(\Omega)}}{||{\Phi_h}(t)||_{L^2(\Omega)}},
\end{equation}
where $\Phi_h$ and $\Phi_r$ are the FOM approximation of a given
field (i.e., $\v_h$, $\u_h$, $q_h$ or $\qbar_h$) and the corresponding ROM 
approximation (i.e., $\v_r$, $\u_r$, $q_r$ or $\qbar_r$), respectively. 
Fig.~\ref{fig:err_t} shows error \eqref{eq:error1} for the two velocity and pressure fields over time 
for the three different ROM techniques under investigation: PPE, $\text{SUP}_1$ and $\text{SUP}_2$. 
We also present the errors related to the ROM computations with no stabilization technique, referred to 
as NOS.
As one would expect, Fig.~\ref{fig:err_t} shows that the model with no stabilization 
is completely unreliable. 
From Fig.~\ref{fig:err_t}, we see that the $\text{SUP}_1$-ROM is a big improvement
over NOS-ROM. However, while the errors for the velocity fields are acceptable,
 the pressure errors remain large and far above the values obtained with NSE in \cite{Stabile2018}. 
 This shows that $\text{SUP}_1$-ROM is not a reliable stabilization of the ROM for the EF algorithm. 
 The $\text{SUP}_2$-ROM produces much better results in terms of the pressure fields, 
 with errors for the velocity fields that are comparable with those given by the $\text{SUP}_1$-ROM. 
 We speculate that the better performance of the $\text{SUP}_2$ model  (which adds 
 the supremizer solutions related to \emph{both} pressure fields to the evolve and filtered velocity spaces) 
 with respect to the $\text{SUP}_1$ model could be due to the strong coupling between the evolve velocity 
 $\v$ and the filter velocity $\u$, 
However, the stability of the inf-sup ROM formulation for the EF algorithm needs to be investigate in more depth
and will be object of future work. 
Finally, from  Fig.~\ref{fig:err_t} we observe that the PPE-ROM provides the lowest errors
for the velocity fields, while the pressure errors are comparable (for $q$)
or worse (for $\overline{q}$) than the errors given by the $\text{SUP}_2$-ROM.
Our findings for the EF algorithm are in agreement with what observed in 
\cite{Stabile2018} for NSE. Indeed, therein it is shown that the PPE model
produces better results for the velocity field but worse results for the pressure field
when compared to the supremizer enrichment model.  
For insights on the behavior of the errors at the first and last time steps of the simulation, 
we refer to \cite{Girfoglio2020}. 

We report minimum, average, and maximum relative errors for PPE and $\text{SUP}_2$ in Table \ref{tab:errors_t}. 
The average errors for $\u$ and $q$ are comparable to the values obtained in \cite{Stabile2018}.

\begin{table}[htb!]
\centering
\footnotesize
\begin{tabular}{lcccccccc}
\multicolumn{2}{c}{} \\
\cline{1-9}
 & $\v$ (PPE) & $\v$ ($\text{SUP}_2$) & $\u$ (PPE) & $\u$ ($\text{SUP}_2$) & $q$ (PPE) & $q$ ($\text{SUP}_2$) & $\overline{q}$ (PPE) & $\overline{q}$ ($\text{SUP}_2$)  \\
\hline
Maximum $E_\Phi$ & 9.1e-1  & 9.1e-1 & 9.2e-1 & 9.2e-1 & 3.6e-1 & 5.8e-1 & 5.4e-1 & 9.4e-1 \\
Average $E_\Phi$ & 2.3e-2 & 2.6e-2 & 2.4e-2 & 2.6e-2 & 1.4e-1 & 1.7e-1 & 1.3e-1 & 6e-2\\
Minimum $E_\Phi$ & 7.8e-4 & 6.3e-3 & 7.8e-4 & 6.3e-3 & 2.7e-2 & 4.4e-2 & 1.1e-1 & 1.7e-2 \\
\hline
\end{tabular}
\caption{Maximum, average, and minimum relative errors for the velocity and pressure fields for PPE-ROM and $\text{SUP}_2$-ROM.}
\label{tab:errors_t}
\end{table}

\begin{figure}[htb!]
\centering
       \begin{overpic}[width=0.3\textwidth]{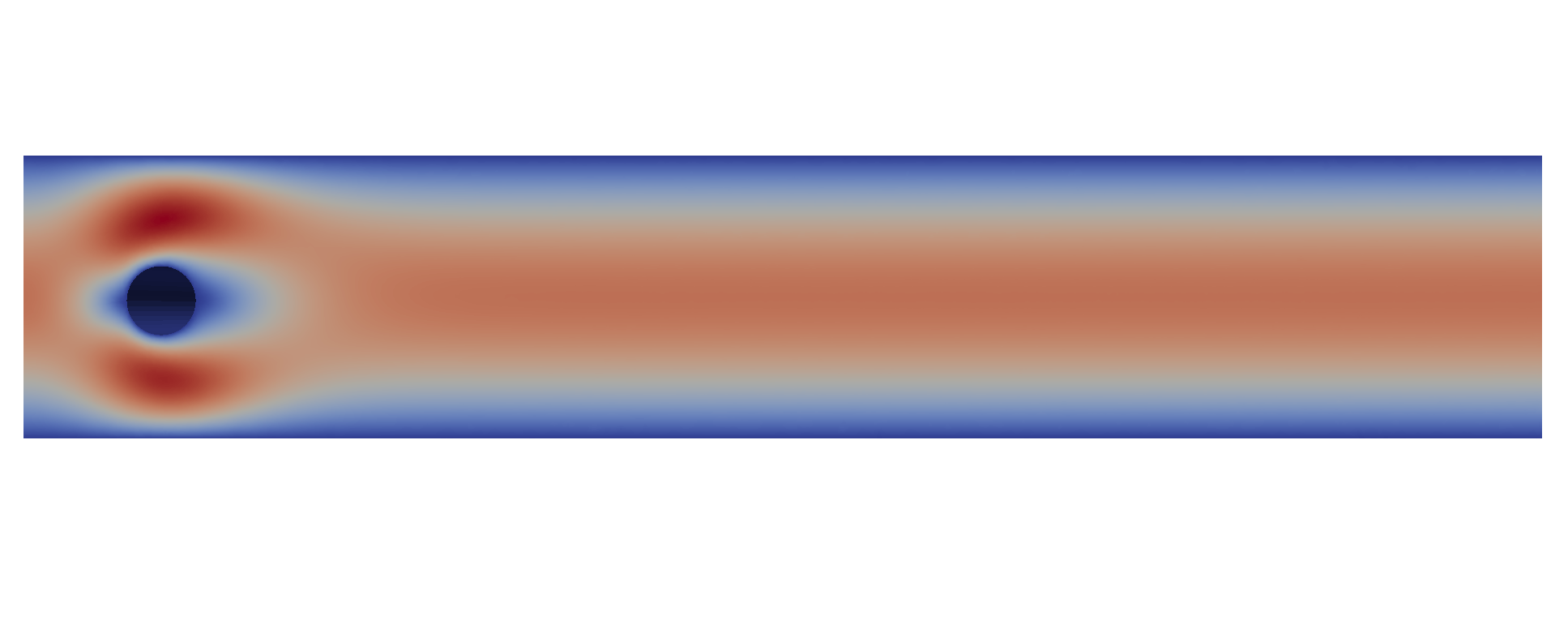}
        \put(40,35){$\v_{FOM}$}
      \end{overpic}
       \begin{overpic}[width=0.3\textwidth]{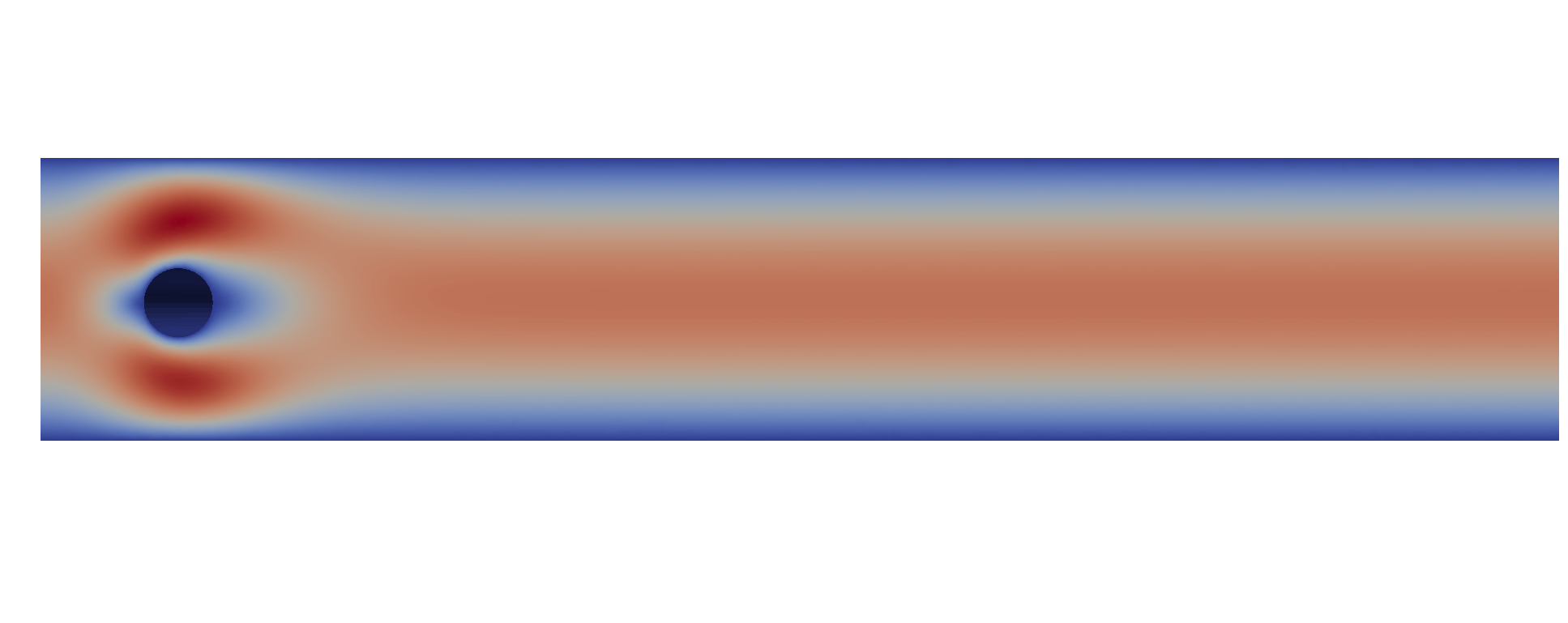}
        \put(40,35){$\v_{PPE}$}
      \end{overpic} 
       \begin{overpic}[width=0.3\textwidth]{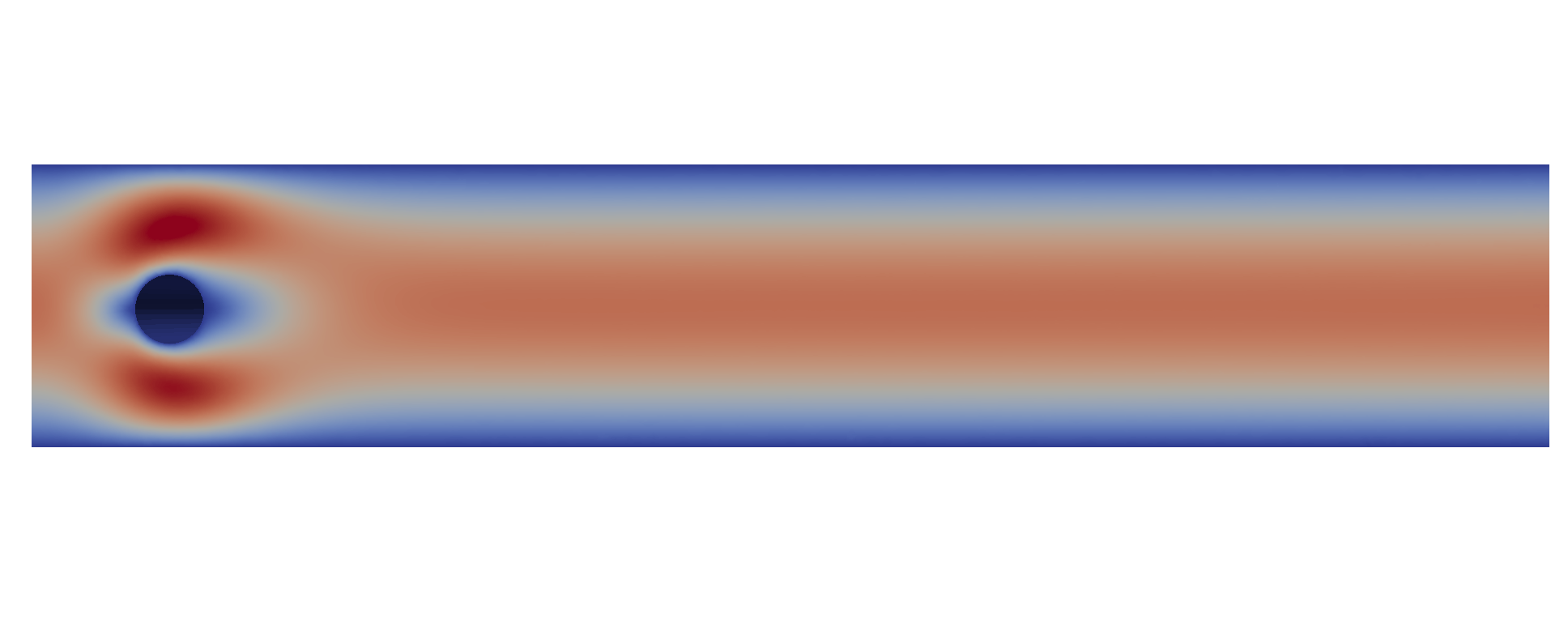}
        \put(40,35){$\v_{SUP_2}$}
      \end{overpic}
\begin{overpic}[width=0.055\textwidth]{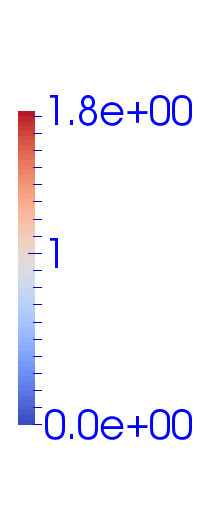}
      \end{overpic}
       \begin{overpic}[width=0.3\textwidth]{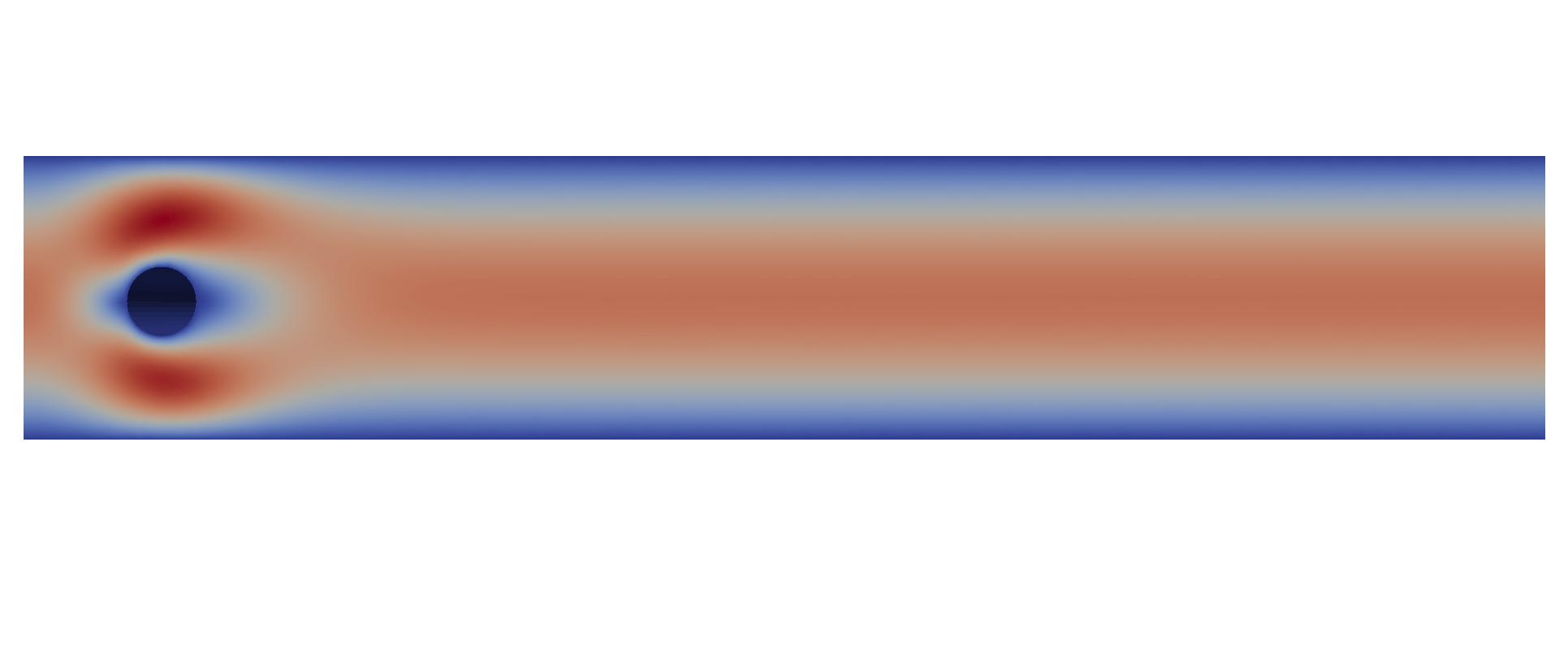}
        \put(40,35){$\u_{FOM}$}
      \end{overpic}
       \begin{overpic}[width=0.3\textwidth]{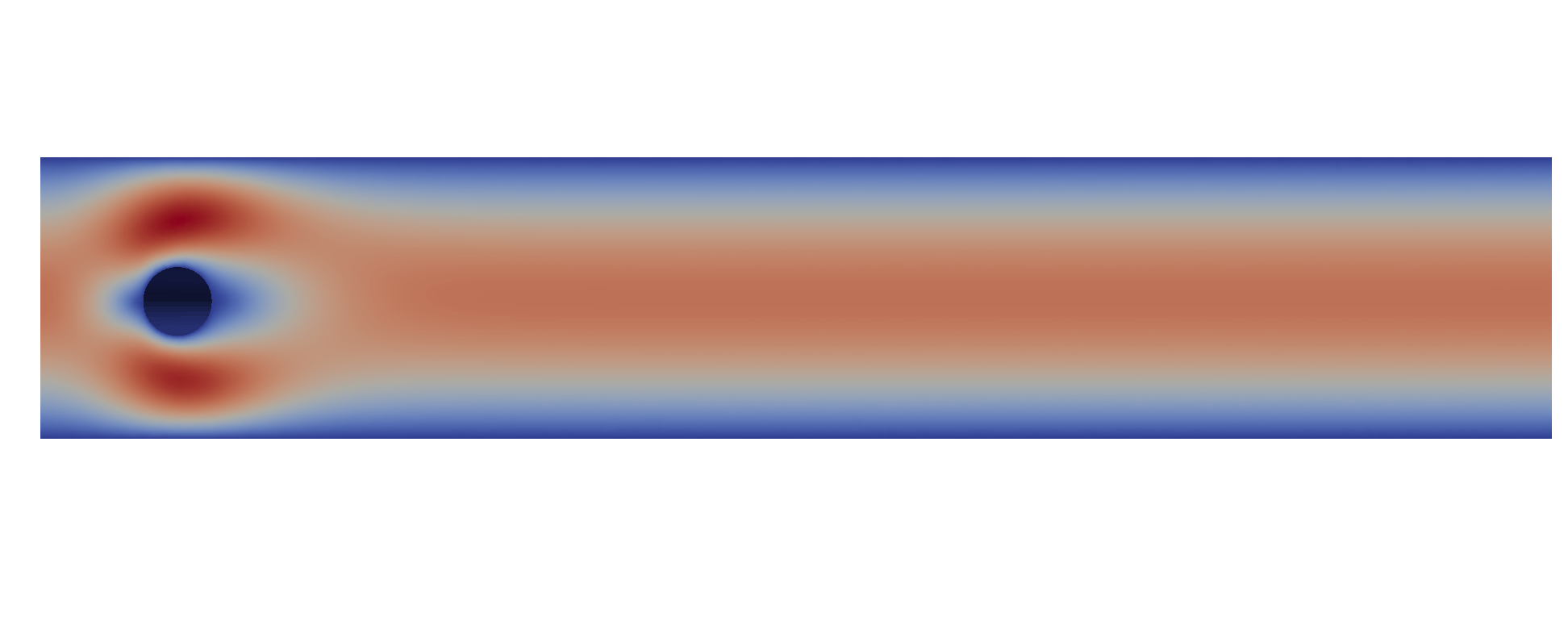}
        \put(40,35){$\u_{PPE}$}
      \end{overpic}
       \begin{overpic}[width=0.3\textwidth]{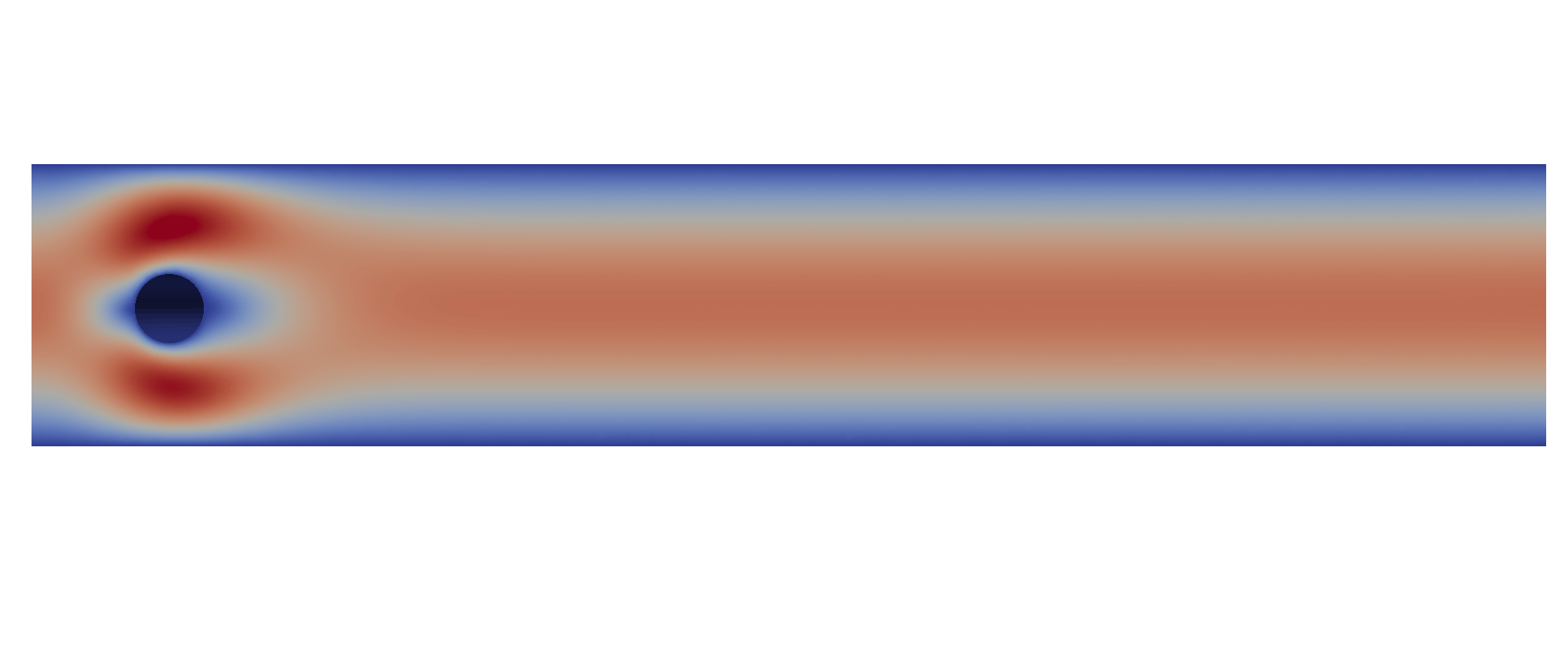}
        \put(40,35){$\u_{SUP_2}$}
      \end{overpic}
\begin{overpic}[width=0.055\textwidth]{img/legendU.png}
      \end{overpic}
       \begin{overpic}[width=0.3\textwidth]{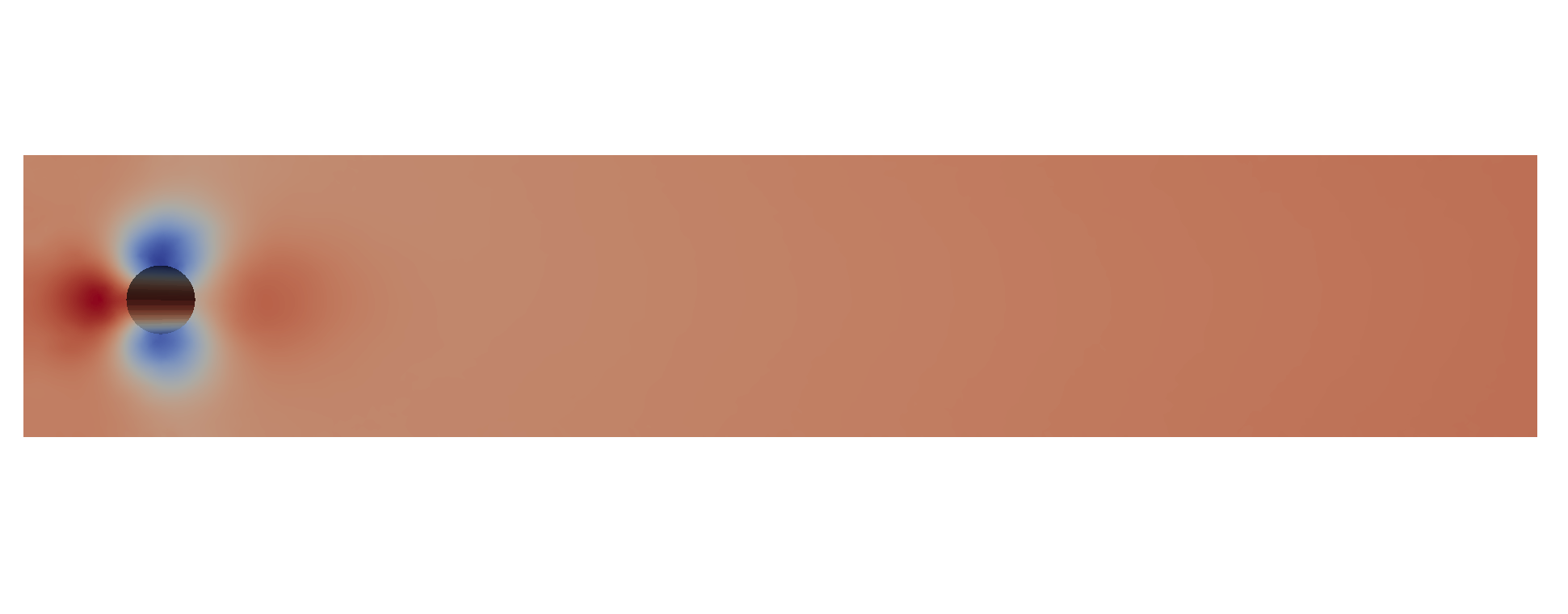}
        \put(40,35){$q_{FOM}$}
      \end{overpic}
       \begin{overpic}[width=0.3\textwidth]{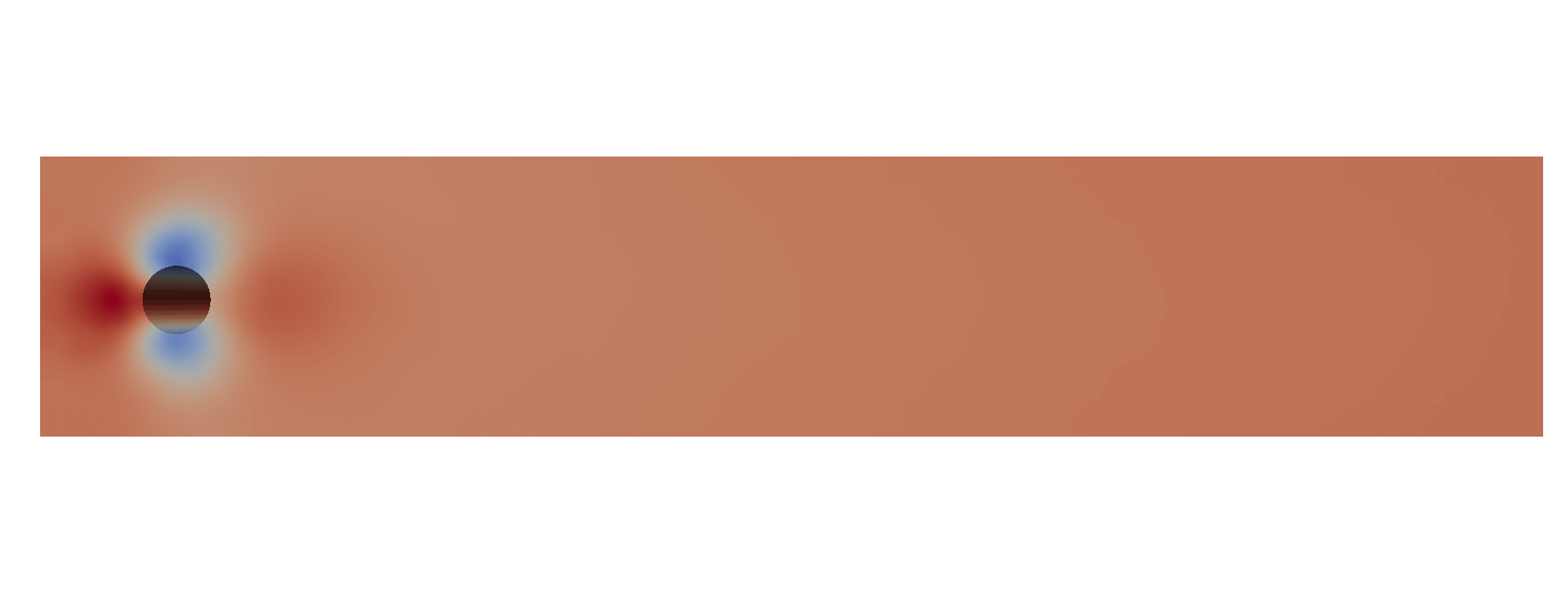}
        \put(40,35){$q_{PPE}$}
      \end{overpic}
 \begin{overpic}[width=0.3\textwidth]{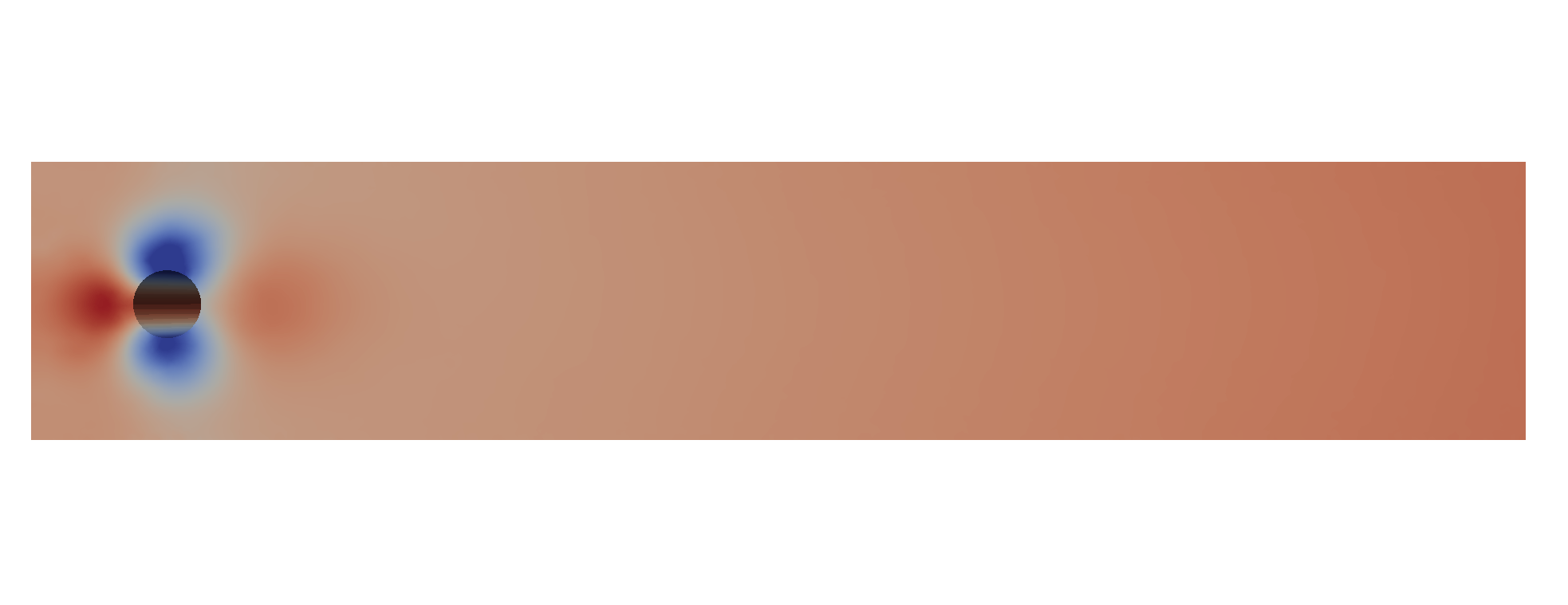}
        \put(40,35){$q_{SUP_2}$}
      \end{overpic}
\begin{overpic}[width=0.055\textwidth]{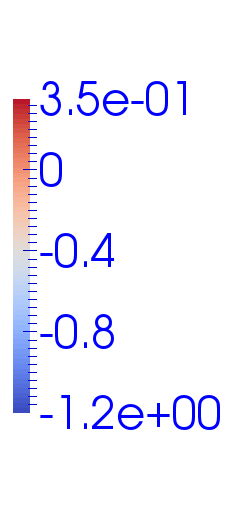}
      \end{overpic}
       \begin{overpic}[width=0.3\textwidth]{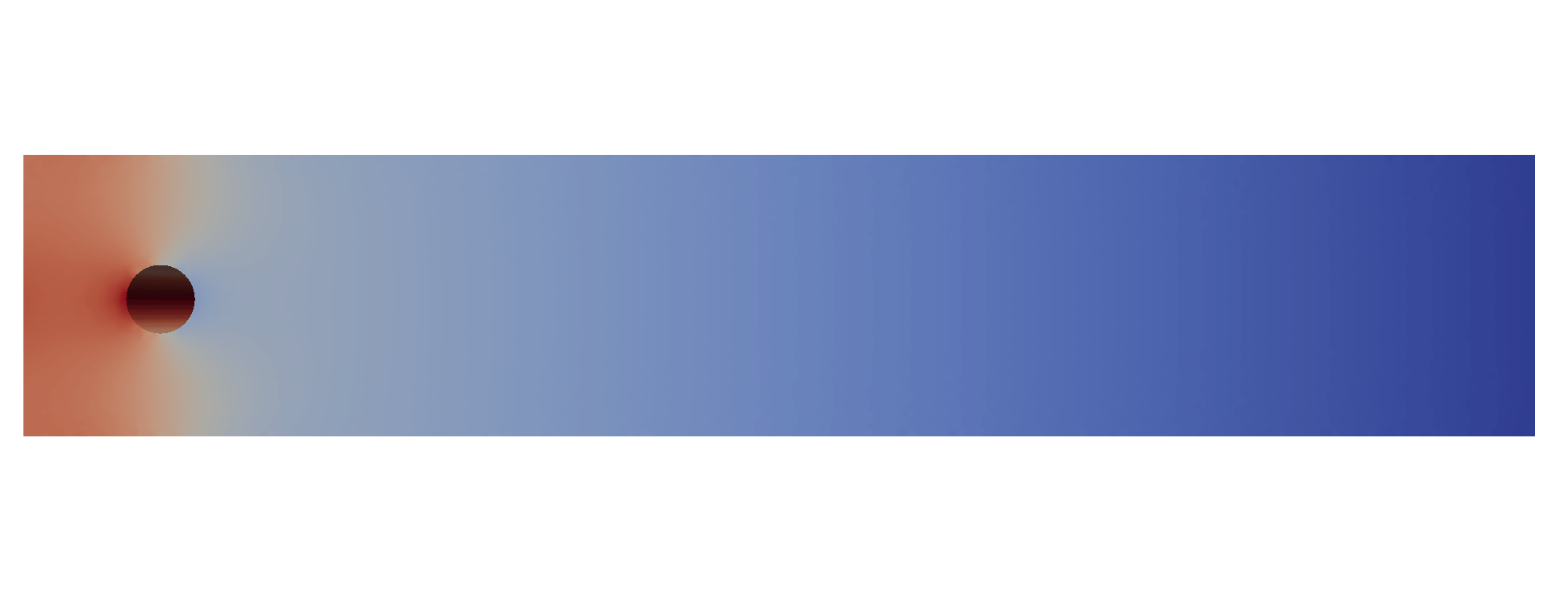}
        \put(40,35){$\overline{q}_{FOM}$}
      \end{overpic}
       \begin{overpic}[width=0.3\textwidth]{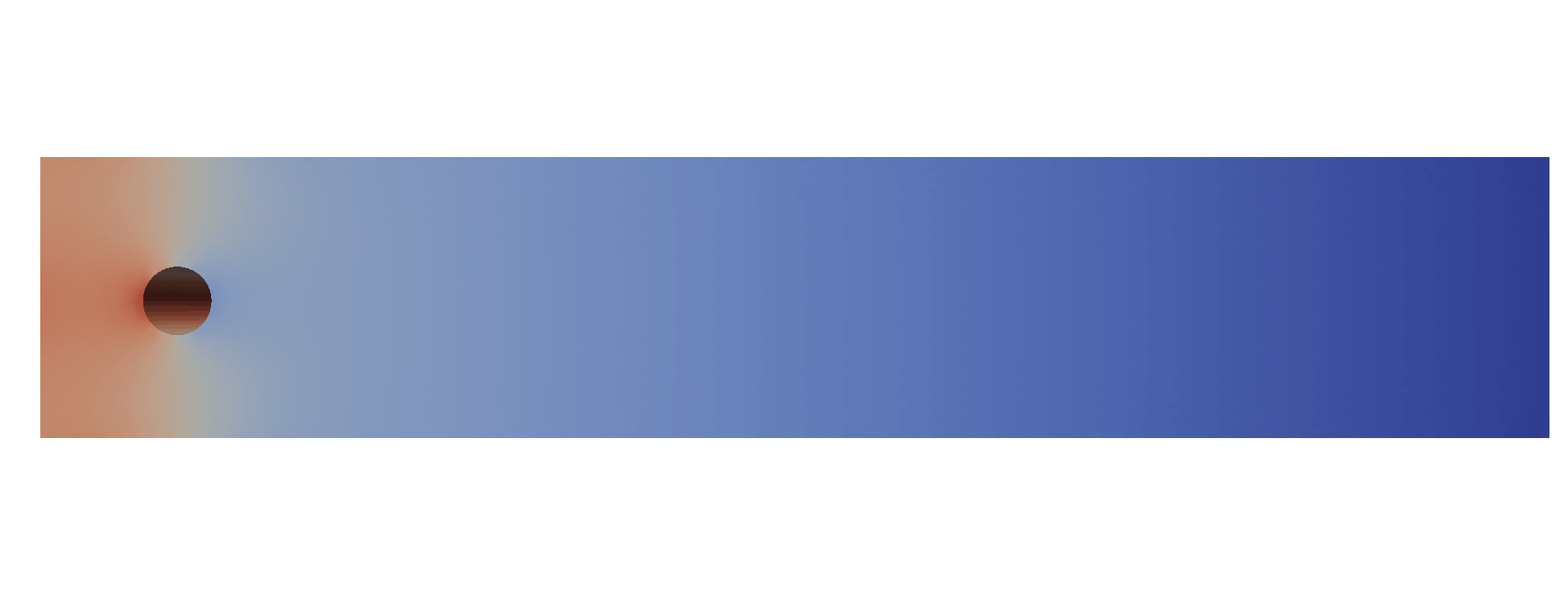}
        \put(40,35){$\overline{q}_{PPE}$}
      \end{overpic}
       \begin{overpic}[width=0.3\textwidth]{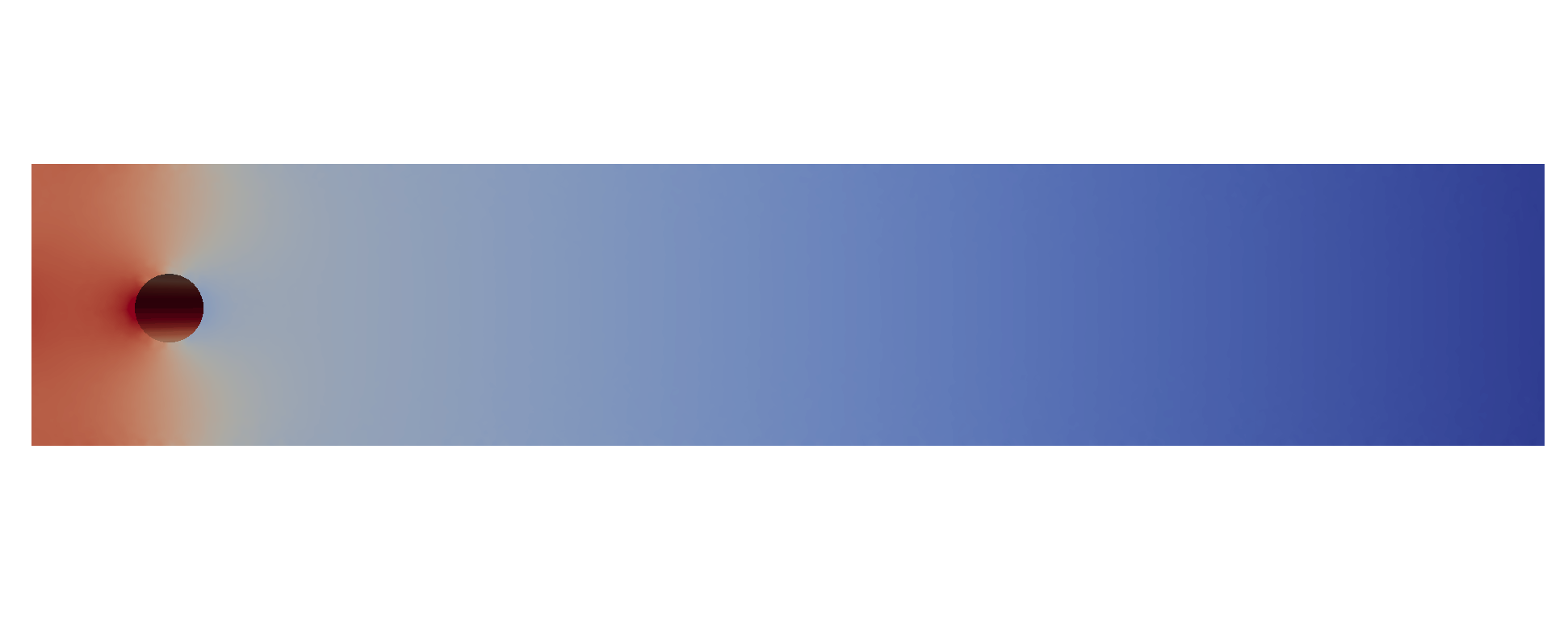}
        \put(40,35){$\overline{q}_{SUP_2}$}
      \end{overpic}
\begin{overpic}[width=0.05\textwidth]{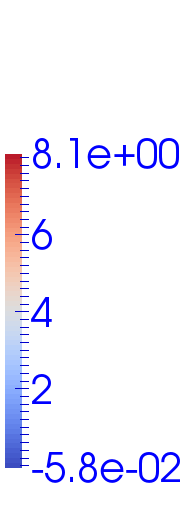}
      \end{overpic}
\caption{
Velocity fields $\u$ (first row) and $\v$ (second row) and pressure fields $q$ (third row)
and $\overline{q}$ (fourth row) at time $t = 5$ computed by FOM (left), PPE-ROM (center), and 
$\text{SUP}_2$-ROM (right).
}
\label{fig:FOMROM}
\end{figure}

For a visual comparison, we report in Fig.~\ref{fig:FOMROM} velocity and pressure fields at $t = 5$ computed by
FOM, PPE-ROM, and $\text{SUP}_2$-ROM. We observe that both PPE-ROM and $\text{SUP}_2$-ROM 
can capture well the main flow features.

\begin{figure}[htb!]
\centering
 \begin{overpic}[width=0.45\textwidth]{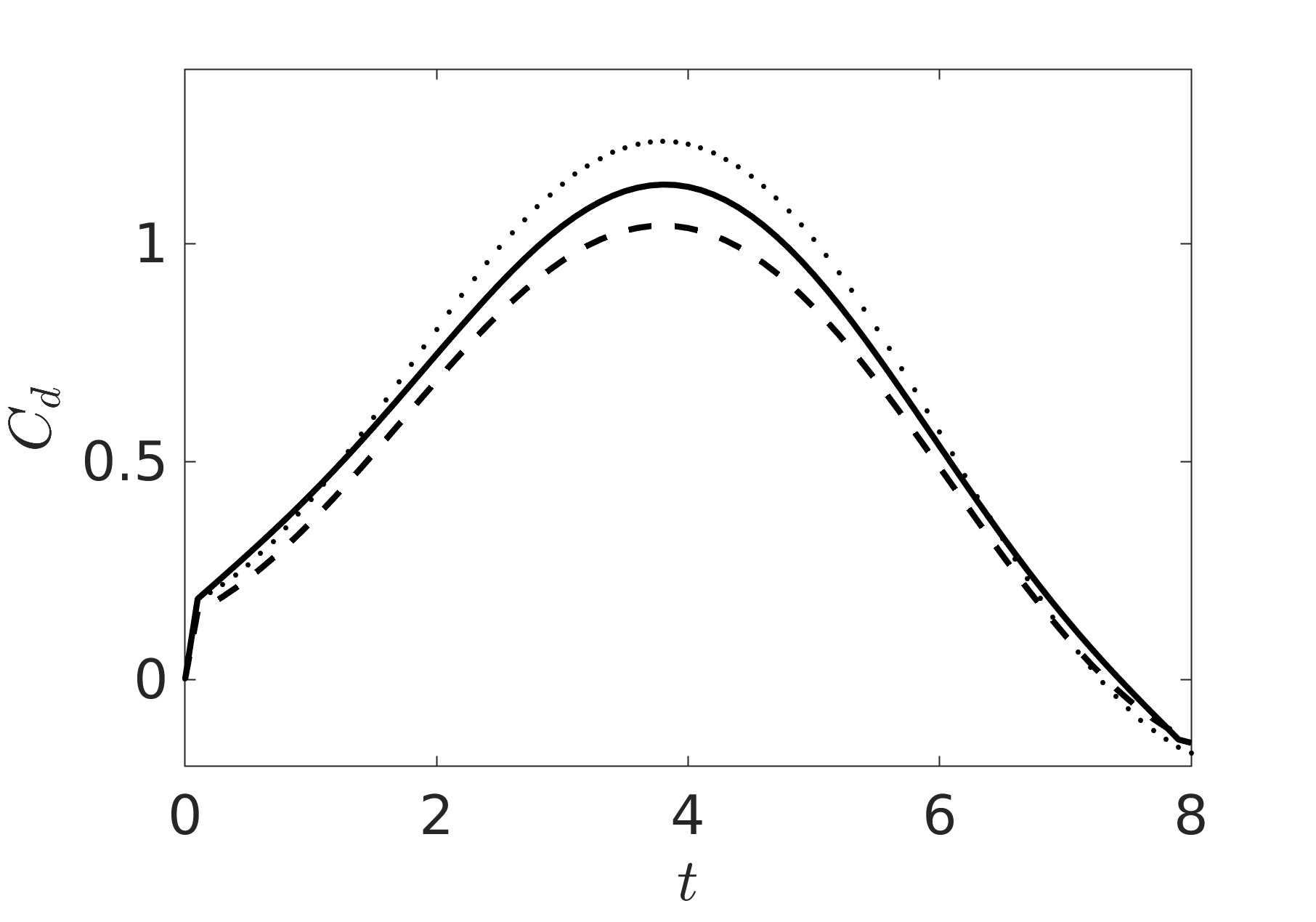}
      \end{overpic}
       \begin{overpic}[width=0.45\textwidth]{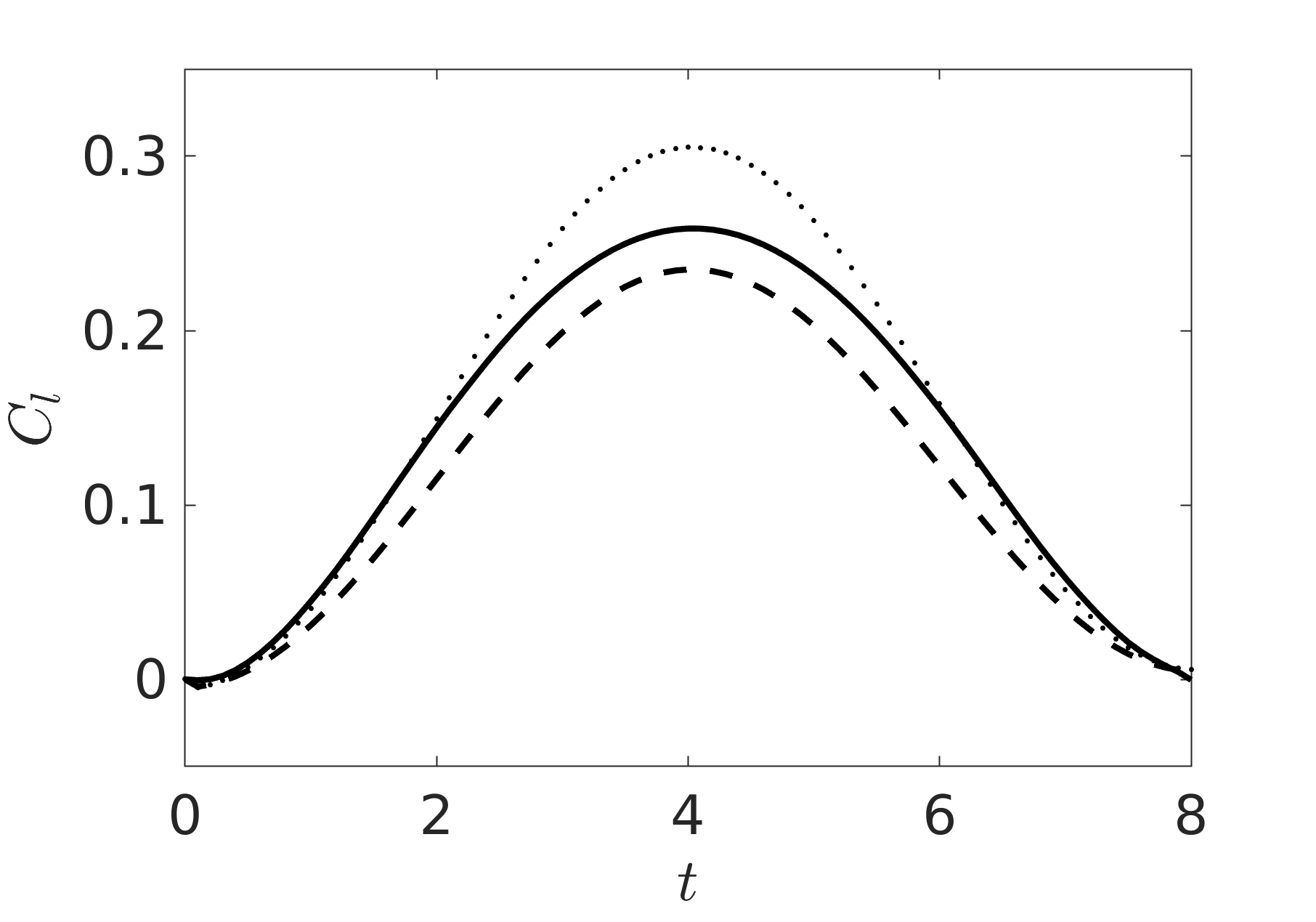}
      \end{overpic}\\
\caption{Aerodynamic coefficients $C_d$ (left) and $C_l$ (right) computed by FOM (solid line), $\text{SUP}_2$-ROM (dotted line), 
and PPE-ROM (dashed line). 
}
\label{fig:coeff_t}
\end{figure}

For a quantitative comparison of FOM, PPE-ROM, and $\text{SUP}_2$-ROM, we consider
the quantities of interest for this benchmark, i.e.~the drag and lift coefficients \cite{John2004,turek1996}:
\begin{align}\label{eq:cd_cl}
c_d(t) = \dfrac{2}{\rho L_{r} {U}^2_{r}} \int_S \left(\left(2 \mu \nabla \u - q\boldsymbol{I}\right)
\cdot \boldsymbol{n}\right) \cdot \boldsymbol{t}~dS, \quad
c_l(t) = \dfrac{2}{\rho L_{r} {U}^2_{r}} \int_S \left(\left(2 \mu \nabla \u - q\boldsymbol{I}\right)
\cdot \boldsymbol{n}\right) \cdot \boldsymbol{n}~dS,
\end{align}
where $U_{r}= 1$ is the maximum velocity at the inlet/outlet, $L_r = 0.1$ is the cylinder diameter, 
$S$ is the cylinder surface, and $\boldsymbol{t}$ and $\boldsymbol{n}$ are the tangential and normal unit vectors
to the cylinder, respectively. 
See Fig.~\ref{fig:coeff_t} for the coefficients in \eqref{eq:cd_cl} computed by the three approaches.
We observe that the amplitude of both coefficients is slightly underestimated (resp., overestimated) 
by the PPE-ROM (resp., $\text{SUP}_2$-ROM) over the entire time interval.
The ROM reconstruction of the lift coefficient appears to be more critical, especially for the $\text{SUP}_2$-ROM and 
around the center of the time interval. 

For a further quantitative assessment of the reconstruction of the coefficients in \eqref{eq:cd_cl}, 
we computed the following errors
\begin{align}\label{eq:error_coeff}
E_{c_d} = \dfrac{||c_d(t)^{FOM}- c_d(t)^{ROM}||_{L^2(0,8)}}{||c_d(t)^{FOM}||_{L^2(0,8)}}, \quad
E_{c_l} = \dfrac{||c_l(t)^{FOM}- c_l(t)^{ROM}||_{L^2(0,8)}}{||c_l(t)^{FOM}||_{L^2(0,8)}}.
\end{align}
We notice that errors defined in \eqref{eq:error_coeff} are different from the ones 
in \cite{Girfoglio2020}, which are related to the maximum values only.
Table \ref{tab:coeffs_t} reports the errors \eqref{eq:error_coeff} for PPE-ROM and $\text{SUP}_2$-ROM. We see that the two ROM strategies provide comparable results: slightly lower than 9\% relative error for the
 the drag coefficient and around 14\% relative error for the lift coefficient. 


\begin{table}[htp!]
\centering
\begin{tabular}{|c|c|c|}
\hline
& $E_{c_d}$ & $E_{c_l}$ \\
\hline
PPE-ROM & 0.088 & 0.144  \\
\hline
$\text{SUP}_2$-ROM & 0.086 & 0.139 \\
\hline
\end{tabular}
\caption{Relative errors \eqref{eq:error_coeff} for lift and drag coefficients for PPE-ROM and $\text{SUP}_2$-ROM.}
\label{tab:coeffs_t}
\end{table}

Finally, we provide some information on the computational cost. The total CPU time required by a FOM simulation is about $1900$ s. 
The solution of the PPE-ROM algebraic systems \eqref{eq:reduced_1}, \eqref{eq:reduced_q}, \eqref{eq:reduced2_1}, and \eqref{eq:reduced_q_bar} takes $1.61$ s, while the solution of the $\text{SUP}_2$-ROM systems \eqref{eq:reduced_1}, \eqref{eq:reduced_2}, \eqref{eq:reduced2_1}, and \eqref{eq:reduced2_2} takes $2.56$ s. 
So, the $\text{SUP}_2$-ROM is less efficient than the PPE-ROM. 
The additional cost of the $\text{SUP}_2$-ROM comes from the 
supremizer modes, which increases the size of the reduced dynamical system. However, it is possible to assert that both ROMs allow to obtain a considerable speed-up.


\section{Conclusions and future perspectives}\label{sec:conclusions}

We presented a POD-Galerkin based reduced order method for a Leray model implemented through the
Evolve-Filter (EF) algorithm. 
Unlike the large majority of the works on Leray-type models, we choose a 
Finite Volume method for the space discretizaion because of its computational efficiency.  
The novelty of this work is the investigation and comparison of 
two techniques for the stabilization of the pressure fields: (i) Poisson Pressure equation, and 
(ii) supremizer enrichment. We showed that the standard supremizer enrichment, 
which works well for the Navier-Stokes equations with no filter, 
needs to be modified in order to obtain stable and accurate solutions with the EF algorithm. 
The modification consists in adding to the evolve and filter velocity spaces the supremizer solutions
related to both evolve and filter pressure fields. 
We assessed our ROM through the classical 2D flow past a cylinder benchmark.
We found that our ROM with both Poisson Pressure equation and modified 
supremizer enrichment captures the flow features with an accuracy comparable to ROMs
applied to the Navier-Stokes equations with no filter \cite{Stabile2018}. Moreover, 
we quantified the relative error of the drag and lift coefficients computed by ROM and FOM
and found that both stabilization approaches produce comparable errors. 

In the future, we would like to investigate in more depth the inf-sup stability of a ROM formulation 
with supremizer enrichment for the EF algorithm. In addition, we would like to extend to the 
EF algorithm other efficient stabilization techniques, such as the one proposed 
in \cite{StabileZancanaroRozza2020}. 

\section*{Disclosure statement}
No potential conflict of interest was reported by the author(s).

\section*{Funding}\label{sec:acknowledgements}
This research was funded by the European Research Council
Executive Agency by the Consolidator Grant project AROMA-CFD "Advanced
Reduced Order Methods with Applications in Computational Fluid Dynamics"
- GA 681447, H2020-ERC CoG 2015 AROMA-CFD, PI G. Rozza, INdAM-
GNCS 2019-2020 projects, the US National Science Foundation through grant DMS-1620384 and DMS-195353.

\bibliographystyle{plain}
\bibliography{FSI.bib}

\begin{thebibliography}{10}

\bibitem{akhtar2009stability}
I.~Akhtar, A.~H. Nayfeh, and C.~J. Ribbens.
\newblock On the stability and extension of reduced-order galerkin models in
  incompressible flows.
\newblock {\em Theoretical and Computational Fluid Dynamics}, 23(3):213--237,
  2009.

\bibitem{Aubry1988}
N.~Aubry, P.~Holmes, J.~L. Lumley, and E.~Stone.
\newblock The dynamics of coherent structures in the wall region of a turbulent
  boundary layer.
\newblock {\em Journal of Fluid Mechanics}, 192(-1):115, jul 1988.

\bibitem{Bader2016}
E.~Bader, M.K\"{a}rcher, M.~A. Grepl, and K.~Veroy.
\newblock {Certified Reduced Basis Methods for Parametrized Distributed
  Elliptic Optimal Control Problems with Control Constraints}.
\newblock {\em {SIAM} Journal on Scientific Computing}, 38(6):A3921--A3946, jan
  2016.

\bibitem{Ballarin2014}
F.~Ballarin, A.~Manzoni, A.~Quarteroni, and G.~Rozza.
\newblock {Supremizer stabilization of {POD}-Galerkin approximation of
  parametrized steady incompressible Navier-Stokes equations}.
\newblock {\em International Journal for Numerical Methods in Engineering},
  102(5):1136--1161, nov 2014.

\bibitem{Benner2015}
P.~Benner, S.~Gugercin, and K.~Willcox.
\newblock {A Survey of Projection-Based Model Reduction Methods for Parametric
  Dynamical Systems}.
\newblock {\em {SIAM} Review}, 57(4):483--531, jan 2015.

\bibitem{bennerParSys}
P.~Benner, M.~Ohlberger, A.~Patera, and K.~Rozza, G.and~Urban.
\newblock {\em {Model Reduction of Parametrized Systems}}, volume 1st ed. 2017
  of {\em MS\&A series}.
\newblock Springer, 2017.

\bibitem{ModelOrderReduction}
P.~Benner, W.~Schilders, S.~Grivet-Talocia, A.~Quarteroni, G.~Rozza, and L.~M.
  Silveira.
\newblock {\em Model Order Reduction}.
\newblock De Gruyter, Berlin, Boston, 2020.

\bibitem{Bergmann2009Enablers}
M.~Bergmann, C.-H. Bruneau, and A.~Iollo.
\newblock Enablers for robust {POD} models.
\newblock {\em Journal of Computational Physics}, 228(2):516--538, feb 2009.

\bibitem{BQV}
L.~Bertagna, A.~Quaini, and A.~Veneziani.
\newblock {Deconvolution-based nonlinear filtering for incompressible flows at
  moderately large Reynolds numbers}.
\newblock {\em International Journal for Numerical Methods in Fluids},
  81(8):463--488, 2016.

\bibitem{boffi_mixed}
D.~Boffi, F.~Brezzi, and M.~Fortin.
\newblock {\em {Mixed Finite Element Methods and Applications}}.
\newblock Springer-Verlag Berlin Heidelberg, first edition, 2013.

\bibitem{abigail_CMAME}
A.L. Bowers and L.G. Rebholz.
\newblock Numerical study of a regularization model for incompressible flow
  with deconvolution-based adaptive nonlinear filtering.
\newblock {\em Computer Methods in Applied Mechanics and Engineering},
  258:1--12, 2013.

\bibitem{Boyd1998283}
J.~P. Boyd.
\newblock Two comments on filtering (artificial viscosity) for {Chebyshev} and
  {Legendre} spectral and spectral element methods: Preserving boundary
  conditions and interpretation of the filter as a diffusion.
\newblock {\em Journal of Computational Physics}, 143(1):283 -- 288, 1998.

\bibitem{BREZZI199027}
F.~Brezzi and K.-J. Bathe.
\newblock {A discourse on the stability conditions for mixed finite element
  formulations}.
\newblock {\em Computer Methods in Applied Mechanics and Engineering}, 82(1):27
  -- 57, 1990.
\newblock Proceedings of the Workshop on Reliability in Computational
  Mechanics.

\bibitem{Carlberg2013623}
K.~Carlberg, C.~Farhat, J.~Cortial, and D.~Amsallem.
\newblock {The GNAT method for nonlinear model reduction: Effective
  implementation and application to computational fluid dynamics and turbulent
  flows}.
\newblock {\em Journal of Computational Physics}, 242:623 -- 647, 2013.

\bibitem{ChinestaEnc2017}
F.~Chinesta, A.~Huerta, G.~Rozza, and K.~Willcox.
\newblock {Model Order Reduction}.
\newblock {\em Encyclopedia of Computational Mechanics, Elsevier Editor}, 2016.

\bibitem{Chinesta2011}
F.~Chinesta, P.~Ladeveze, and E.~Cueto.
\newblock {A short review on Model Order Reduction based on Proper Generalized
  Decomposition}.
\newblock {\em Archives of Computational Methods in Engineering}, 18(4):395,
  2011.

\bibitem{Dumon20111387}
A.~Dumon, C.~Allery, and A.~Ammar.
\newblock {Proper General Decomposition (PGD) for the resolution of
  Navier-Stokes equations}.
\newblock {\em Journal of Computational Physics}, 230(4):1387--1407, 2011.

\bibitem{Dunca2005}
A.~Dunca and Y.~Epshteyn.
\newblock On the {Stolz-Adams} deconvolution model for the large-eddy
  simulation of turbulent flows.
\newblock {\em SIAM Journal on Mathematical Analysis}, 37(6):1890--1902, 2005.

\bibitem{Fischer2001265}
P.~Fischer and J.~Mullen.
\newblock Filter-based stabilization of spectral element methods.
\newblock {\em Comptes Rendus de l'Academie des Sciences - Series I -
  Mathematics}, 332(3):265 -- 270, 2001.

\bibitem{Gerner2011}
A.-L Gerner and K.~Veroy.
\newblock Certified reduced basis methods for parametrized saddle point
  problems.
\newblock {\em SIAM Journal on Scientific Computing}, 34, 01 2011.

\bibitem{Girfoglio2021}
M.~Girfoglio, A.~Quaini, and G.Rozza.
\newblock Fluid-structure interaction simulations with a les filtering approach
  in solids4foam.
\newblock {\em accepted for the publication in Communications in Applied and
  Industrial Mathematics}, 2021.
\newblock https://arxiv.org/abs/2102.08011.

\bibitem{Girfoglio2019}
M.~Girfoglio, A.~Quaini, and G.~Rozza.
\newblock {A Finite Volume approximation of the Navier-Stokes equations with
  nonlinear filtering stabilization}.
\newblock {\em Computers \& Fluids}, 187:27--45, 2019.

\bibitem{Girfoglio2020}
M.~Girfoglio, A.~Quaini, and G.~Rozza.
\newblock {A POD-Galerkin reduced order model for a LES filtering approach}.
\newblock {\em Journal of Computational Physics}, 436:110260, 2021.

\bibitem{Gunzburger2019}
M.~Gunzburger, T.~Iliescu, M.~Mohebujjaman, and M.~Schneier.
\newblock An evolve-filter-relax stabilized reduced order stochastic
  collocation method for the time-dependent navier-stokes equations.
\newblock {\em SIAM/ASA Journal on Uncertainty Quantification}, 7:1162--1184,
  01 2019.

\bibitem{Gunzburger2019b}
M.~Gunzburger, T.~Iliescu, and M.~Schneier.
\newblock A {Leray} regularized ensemble-proper orthogonal decomposition method
  for parameterized convection-dominated flows.
\newblock {\em IMA Journal of Numerical Analysis}, 40(2):886--913, 2019.

\bibitem{hesthaven2015certified}
J.~S. Hesthaven, G.~Rozza, and B.~Stamm.
\newblock {\em {Certified Reduced Basis Methods for Parametrized Partial
  Differential Equations}}.
\newblock Springer International Publishing, 2016.

\bibitem{PISO}
R.~I. Issa.
\newblock Solution of the implicitly discretised fluid flow equations by
  operator-splitting.
\newblock {\em Journal of Computational Physics}, 62(1):40--65, 1986.

\bibitem{John2004}
V.~John.
\newblock {Reference values for drag and lift of a two dimensional
  time-dependent flow around a cylinder}.
\newblock {\em International Journal for Numerical Methods in Fluids},
  44:777--788, 2004.

\bibitem{JOHNSTON2004221}
H.~Johnston and J.-G. Liu.
\newblock {Accurate, stable and efficient Navier-Stokes solvers based on
  explicit treatment of the pressure term}.
\newblock {\em Journal of Computational Physics}, 199(1):221 -- 259, 2004.

\bibitem{Kalashnikova_ROMcomprohtua}
I.~Kalashnikova and M.~F. Barone.
\newblock {On the stability and convergence of a Galerkin reduced order model
  (ROM) of compressible flow with solid wall and far-field boundary treatment}.
\newblock {\em International Journal for Numerical Methods in Engineering},
  83(10):1345--1375, 2010.

\bibitem{Kunisch2002492}
K.~Kunisch and S.~Volkwein.
\newblock {Galerkin proper orthogonal decomposition methods for a general
  equation in fluid dynamics}.
\newblock {\em SIAM Journal on Numerical Analysis}, 40(2):492--515, 2002.

\bibitem{Lax1960}
P.D. Lax and B.~Wendroff.
\newblock System of conservation laws.
\newblock {\em Communications on Pure and Applied Mathematics}, 13:217--237,
  1960.

\bibitem{layton_CMAME}
W.~Layton, L.G. Rebholz, and C.~Trenchea.
\newblock Modular nonlinear filter stabilization of methods for higher
  {Reynolds} numbers flow.
\newblock {\em Journal of Mathematical Fluid Mechanics}, 14:325--354, 2012.

\bibitem{Leray34}
J.~Leray.
\newblock Essai sur le mouvement d'un fluide visqueux emplissant l'espace.
\newblock {\em Journal de Math\'ematiques Pures et Appliqu\'ees}, 63:193--248,
  1934.

\bibitem{Li2020}
L.~Li.
\newblock A split-step finite-element method for incompressible navier-stokes
  equations with high-order accuracy up-to the boundary.
\newblock {\em Journal of Computational Physics}, 408(3):213--237, 2020.

\bibitem{Lorenzi2016}
S.~Lorenzi, A.~Cammi, L.~Luzzi, and G.~Rozza.
\newblock {{POD}-Galerkin method for finite volume approximation of
  Navier{\textendash}Stokes and {RANS} equations}.
\newblock {\em Computer Methods in Applied Mechanics and Engineering},
  311:151--179, nov 2016.

\bibitem{Moin1998}
P.~Moin and K.~Mahesh.
\newblock Direct {Numerical} {Simulation}: A tool in turbulence research.
\newblock {\em Annual Review of Fluid Mechanics}, 30(1):539--578, jan 1998.

\bibitem{Orszag1986}
S.~A. Orszag, M.~Israeli, and M.O. Deville.
\newblock {Boundary conditions for incompressible flows}.
\newblock {\em Journal of Scientific Computing}, 1(1):75--111, Mar 1986.

\bibitem{SIMPLE}
S.~V. Patankar and D.~B. Spalding.
\newblock A calculation procedure for heat, mass and momentum transfer in
  three-dimensional parabolic flows.
\newblock {\em International Journal of Heat and Mass Transfer},
  15(10):1787--1806, 1972.

\bibitem{quarteroniRB2016}
A.~Quarteroni, A.~Manzoni, and F.~Negri.
\newblock {\em {Reduced Basis Methods for Partial Differential Equations}}.
\newblock Springer International Publishing, 2016.

\bibitem{quarteroni2007numerical}
A.~Quarteroni, R.~Sacco, and F.~Saleri.
\newblock {\em {Numerical Mathematics}}.
\newblock Springer Verlag, 2007.

\bibitem{Rozza2009}
G.~Rozza.
\newblock {Reduced basis methods for Stokes equations in domains with
  non-affine parameter dependence}.
\newblock {\em Computing and Visualization in Science}, 12(1):23--35, 2009.

\bibitem{Rozza2008}
G.~Rozza, D.~B.~P. Huynh, and A.~T. Patera.
\newblock {Reduced Basis approximation and a Posteriori error estimation for
  affinely parametrized elliptic coercive Partial Differential Equations}.
\newblock {\em Archives of Computational Methods in Engineering}, 15(3):229,
  May 2008.

\bibitem{Rozza2007}
G.~Rozza and K.~Veroy.
\newblock {On the stability of the reduced basis method for Stokes equations in
  parametrized domains}.
\newblock {\em Computer Methods in Applied Mechanics and Engineering},
  196(7):1244 -- 1260, 2007.

\bibitem{RoSta17}
G.~Stabile and G.~Rozza.
\newblock {ITHACA-FV - In real Time Highly Advanced Computational Applications
  for Finite Volumes}.
\newblock Accessed: 2018-01-30.

\bibitem{Stabile2018}
G.~Stabile and G.~Rozza.
\newblock {Finite volume POD-Galerkin stabilised reduced order methods for the
  parametrised incompressible Navier–Stokes equations}.
\newblock {\em Computer \& Fluids}, 173:273--284, 2018.

\bibitem{StabileZancanaroRozza2020}
G.~Stabile, M.~Zancanaro, and G.~Rozza.
\newblock Efficient geometrical parametrization for finite-volume based reduced
  order methods.
\newblock {\em International Journal for Numerical Methods in Engineering},
  121(12):2655--2682, 2020.

\bibitem{Huerta2020}
V.~Tsiolakis, M.~Giacomini, R.~Sevilla, C.~Othmer, and A.~Huerta.
\newblock Parametric solutions of turbulent incompressible flows in openfoam
  via the proper generalised decomposition.
\newblock 06 2020.

\bibitem{turek1996}
S.~Turek and M.~Sch\"afer.
\newblock Benchmark computations of laminar flow around cylinder.
\newblock In E.H. Hirschel, editor, {\em Flow Simulation with High-Performance
  Computers II}, volume~52 of {\em Notes on Numerical Fluid Mechanics}. Vieweg,
  1996.

\bibitem{Doormaal1984}
J.~P. Van~Doormaal and G.~D. Raithby.
\newblock Enhancements of the simple method for predicting incompressible fluid
  flows.
\newblock {\em Numerical Heat Transfer}, 7(2):147--163, 1984.

\bibitem{wang_turb}
Z.~Wang, I.~Akhtar, J.~Borggaard, and T.~Iliescu.
\newblock {Proper orthogonal decomposition closure models for turbulent flows:
  A numerical comparison}.
\newblock {\em Computer Methods in Applied Mechanics and Engineering},
  237–240:10 -- 26, 2012.

\bibitem{Weller1998}
H.~G. Weller, G.~Tabor, H.~Jasak, and C.~Fureby.
\newblock A tensorial approach to computational continuum mechanics using
  object-oriented techniques.
\newblock {\em Computers in physics}, 12(6):620--631, 1998.

\bibitem{Wells2017}
D.~Wells, Z.~Wang, X.~Xie, and T.~Iliescu.
\newblock An evolve-then-filter regularized reduced order model for
  convection-dominated flows.
\newblock {\em International Journal for Numerical Methods in Fluids},
  84:598--615, 2017.

\bibitem{Xie2016}
X.~Xie, D.~Wells, Z.~Wang, and T.~Iliescu.
\newblock Approximate deconvolution reduced order modeling.
\newblock {\em Computer Methods in Applied Mechanics and Engineering},
  313:512--534, 2016.

\bibitem{Xie2018}
X.~Xie, D.~Wells, Z.~Wang, and T.~Iliescu.
\newblock {Numerical analysis of the Leray reduced order model}.
\newblock {\em Journal of Computational and Applied Mathematics}, 328:12--29,
  2018.

\end{thebibliography}





\end{document}